\documentclass[11pt]{article}

\usepackage{amsmath, amssymb, amsfonts, amscd, epsfig, graphicx, color, subfig, hyperref,cite,authblk}

\newcommand{\R}{\mathbb R}

\newcommand{\q}{q}

\newcommand{\C}{\mathbb C}
\newcommand{\G}{\mathbb G}

\newcommand{\M}{\mathcal M}

\newcommand{\N}{{\M}{}_{s_0}}

\newcommand{\rv}{\rho_{0}}
\newcommand{\Gv}{\Gamma_{0}}
\newcommand{\Lv}{L_{0}}

\newcommand{\ro}{\rho}
\newcommand{\Go}{\Gamma}
\newcommand{\Lo}{L}

\newcommand{\RP}{\psi_{z}^{\sigma}}

\newcommand{\vsp}[1]{\textnormal{span}\{#1\}}

\title{A cortical-inspired geometry for contour perception and motion integration}

\author{D. Barbieri\footnote{Institut des Syst\`{e}mes Complexes, Paris \^{I}le-de-France. 57-59 Rue Lhomond, 75005 Paris, France. {\tt davide.barbieri8@gmail.com}}, G. Citti\footnote{Department of Mathematics, University of Bologna. Piazza di Porta San Donato 5, 40126 Bologna, Italy {\tt giovanna.citti@unibo.it}}, G. Cocci\footnote{DEI - Department of Electrical, Electronic, and Information Engineering "Guglielmo Marconi". University of Bologna, Viale del Risorgimento 2, 40136 Bologna, Italy. {\tt giacomo.cocci2@unibo.it}}, A. Sarti\footnote{CAMS - Centre d'Analyse et de Math\'{e}matique Sociales. EHESS, 190-198 Avenue de France, 75244 Paris, France.{\tt alessandro.sarti@ehess.fr}}}

\begin{document}

\maketitle

\begin{abstract}
In this paper we develop a geometrical model of functional architecture for the processing of spatio-temporal visual stimuli. The model arises from the properties of the receptive field linear dynamics of orien\-ta\-tion and speed-selective cells in the visual cortex, that can be embedded in the definition of a geometry where the connectivity between points is driven by the contact structure of a 5D manifold. Then, we compute the stochastic kernels that are the approximations of two Fokker Planck operators associated to the geometry, and implement them as facilitation patterns within a neural population activity model, in order to reproduce some psychophysiological findings about the perception of contours in motion and trajectories of points found in the literature.\\
{\bf keywords}: {Visual cortex \and Lie groups \and Contact geometry \and Galilean group \and Cognitive neuroscience \and Spatio-temporal models}
\end{abstract}

\section{Introduction}\label{sec:intro}
 
The modular structure of the mammalian visual cortex has been discovered in the seventhies by the pioneeristic work of Hubel and Wiesel \cite{Hubel1977,Hubel}. Many families of cells contribute to sampling and coding the stimulus image. Each family is sensible to a specific feature of the image: position, orientation, scale, color, curvature, velocity, stereo.

The main behaviour of simple cells is that of detection of positions and local orientations via linear filtering of the stimulus, and the linear filter associated to a given cell is called its receptive profile. Daugman \cite{Daugman} proposed a criterion of minimal uncertainty for the shape receptive profile, which results in two dimensional spatial Gabor filters, that was later confirmed by Jones and Palmer \cite{Jones1987b} and Ringach \cite{Ringach2002}. Simple cells also possess a temporal behaviour, studied by De Angelis et al. \cite{DeAngelis1993a}, and spatio-temporal receptive profiles could again be interpreted by Cocci et al. \cite{CBS} in terms of minimal uncertainty, resulting in three dimensional Gabor filters. However, not all cell activity in response to stumuli can be justified in terms of linear filters, and nonlinearities can sometimes become relevant to model their behaviours \cite{Graham}.

Cells are spatially organized in such a way that for every point (x,y) of the retinal plane there is an entire set of cells, each one sensitive to a particular instance of the considered feature, giving rise to the so-called hypercolumnar organization. Hypercolumnar organization and neural connectivity between hypercolumns constitute the functional architecture of the visual cortex, that is the cortical structure underlying the processing of visual stimulus.

The mathematical modelling of the functional architecture of the visual cortex in terms of differential geometry was introduced with the seminal works of Koenderink \cite{Koenderink} and Hoffmann \cite{Hoffman66, Hoffman}. While the first author pointed out the differential action of perceptual mechanisms, in particular with respect to jet spaces arising from linear filters, the second author proposed to model the hypercolumnar organization in terms of a fiber bundle structure and pointed out the central role of symmetries in perception expressing them in terms of Lie groups and Lie algebras.

%Around the same time, a different approach was developed by Wilson and Cowan \cite{Wilson1972} and Ermentrout and Cowan \cite{Ermentrout1980}, in terms of dynamics of populations of neural oscillators. Such approach introduced techniques of dynamical systems which make use of nonlinearities from feedback-feedforward mechanisms in the network of connected cells in order to characterize neural morphologies.

%Some years later, Mumford and Shah \cite{MumfordShah} proposed a variational approach to a problem of computer vision that was later recognized to be strictly related to the cited studies. The problem was that of image segmentation, for which they proposed a functional that is singular along the edges of the image.

Problems of perceptions can also be afforded with a purely psychophysical approach. The  study of Field Hayes and Hess \cite{Field1993} introduced the notion of association field, as path of information integration along images that can quantitatively satisfy the assumptions of the Gestalt principle of good continuation. The perceptual role of this mechanism was indeed that of contour integration, that typically occurs along field lines associated to locally coherent directions in images.

Almost simultaneously  Mumford \cite{Mumford} proposed a variational approach to describe smooth edges, in terms of the elastica functional, that could be implemented with stochastic processes defining curves with random curvature at any point. Indeed, they produce probability distributions in the space $\R^2 \times S^1$ of positions and orientations whose probability peaks follow elastica curves. Williams and Jacobs \cite{Williams} could use such stochastic processes to implement a mechanism of stochastic completion, and interpreted the probability kernel they obtained as tensors representing geometric connections on the space of positions and orientations associated to the neural representation of images due to simple cells.

Many of such results dealing with differential geometry were given a unified framework under the new name of \emph{neurogeometry} by Petitot and Tondud \cite{Petitot1999}, who related the association fields of Field Hayess and Hess with the contact geometry introduced by Hoffmann and the elastica of Mumford.

The problem of edge organization in images was then addressed in terms of a stochastic process of the type of Mumford, introducing nonlinearities in order to take into account the role of curvature, by August and Zucker \cite{AugustZucker2000, AugustZucker2003}, while the variational approach of Nitzberg, Mumford and Shiota was extended from edges to level lines of images by Ambriosio and Masnou \cite{AmbrosioMasnou}.

%The up-to-now independent results of Wilson and Cowan on the dynamics of neural oscillators was then related with the variational results of Mumford and Shah by Citti, Manfredini and Sarti, who showed how the phase resulting from the difference equation that regulate the network of neural populations of the first work converges to the flow associated to the functional of the second one.

Then, in \cite{CS}, Citti and Sarti showed how the functional architecture could be described in terms of Lie groups structures. They interpreted the geometric action of receptive profiles as a lifting of level lines into the space $\R^2 \times S^1$ of positions and orientations, and addressed the problem of occlusion with a nonlinear diffusion-concentration process in such a space of liftings. In particular, this approach allows to introduce orientation, instead of depth, as a third dimension for the disentanglement of crossing level lines. Moreover, by making use of the sub-Riemannian structure associated to the Lie symmetries of the Euclidean motion group $\R^2 \ltimes S^1$ they could prove that their geometrical nonlinear diffusion approximates the variational result of Ambrosio and Masnou. In their model, then, contour completion is justified as a propagation in the sub-Riemannian setting, and the integral curves of the vector fields that generate the Lie algebra can be considered as a mathematical representation of the association fields of Field, Hayes and Hess, hence proving the relation between neural mechanisms and image completion. This method was then concretely implemented in \cite{SCS08}.

The problem of boundary completion was also addressed from a slightly different point of view by Zucker \cite{Zucker2006}, who showed the role of Frenet frames.

Exact solutions to the Fokker-Planck equation associated to Mumford stochastic process were provided by Duits and van Almsick \cite{DvA}, and later Duits and Franken \cite{Duits2010a, Duits2010b}  unified such stochastic approach with nonlinear mechanisms of the type of August and Zucker, keeping left invariance with respect to the Lie symmetry of the Euclidean motion group and yet allowing the invertibility of the whole process. Their result was applied to the problem of contour enhancement and contour completion, working on the whole Lie group by means of a representation via suitably defined linear filters. This approach was then extended to different geometric setting by Duits and F\"uhr \cite{Duits2011}, again with applications to the processing of images.

We also mention the works of Hladky and Pauls \cite{Pauls2010} and Boscain et al \cite{BDRS}, where the authors provided technical developments on the structure of the sub-Riemannian paths and minimal surfaces involved in the cited works.

In this paper we propose a mathematical model of cortical functional architecture for the processing of spatio-temporal visual information, that is compatible with both phenomenological experiments and neurophysiological findings. The adopted theoretical framework follows the outlined path of mathematical models of the activity of the visual cortex, and in particular it continues the geometric approach of Citti and Sarti \cite{CS} to a space of higher dimension, since it takes into account time and velocity of stimuli. We will extend the stochastic process of Mumford \cite{Mumford} to this setting, working in a space of liftings arising from the filtering with spatio-temporal receptive profiles, and make use of assumptions like the ones of Ermentrout and Cowan \cite{Ermentrout1980} in order to construct a population dynamics able to provide a new form of association fields adapted to the problem of motion integration and motion completion under occlusion. Moreover, the resulting kernels will be comparable to measured neural activities in the presence of stimuli characterized by their direction of motion.

Our starting point in Section \ref{sec:geometry} is a process of detection resulting from linear filtering with three dimensional Gabor functions with two spatial and one temporal dimensions, which have been proposed as a model of spatio-temporal receptive profiles of primary visual cortex simple cells (see \cite{DeAngelis1993a}, \cite{CBS}). Spatio-temporal Gabor filters extend simple cells Gabor behavior as spatial filters \cite{Daugman,Ringach2002}, that proved its usefulness for the classical task of edge detection, while its role for motion detection was already pointed out in \cite{Petkov2007}. The Gabor transform takes in input a moving image $f(x,t)$, where $x \in \R^2$ are spatial coordinates and $t \in \R$ is time and provides in output  a representation of the signal in the phase space
\begin{displaymath}
f(x,t) \rightarrow F(q,s;p,\nu), 
\end{displaymath}
with $q \in \R^2$ and $s \in \R$ representing spatio-temporal position and $(p,\nu)$ representing, respectively, spatial and temporal frequency.

Since we are mainly interested in spatio-temporal dynamical aspects, we will assign to temporal frequency the meaning of velocity $v$ on the $(x,t)$ coordinates. We will also select the subset of all detected features corresponding to a fixed value of $|p|$. This  will end up to be a 5D manifold, with a contact structure, induced by a normalization of the Liuoville form.  

In Section \ref{sec:diffusion} we will outline that this constraint carries a notion of admissible curves \cite{Geiges} in a deterministic and a stochastic setting, allowing to compute the kernels connecting filters in the 5D manifold. The stochastic processes are completely determined by the described structures of the tangent space of the 5D manifold and  they turn out to be described by the fundamental  solution of a Fokker Planck equation \cite{Oksendal}. In fact the processes contain diffusion in the fiber variables and transport along the remaining generators of admissible tangent directions, in the spirit of \cite{Mumford,Williams}. We will compute these kernels in two limit cases: the motion of a contour at a fixed time instant, and the motion of a point moving in time. The first one reduces to the geometry of contours, with a notion of instantaneous velocity, the second one corresponds to point trajectories in time, and can be related to the subset of the Galilei group \cite{Sorba} on the plane.

In Section \ref{sec:truelife} we will discuss the compatibility of the previously calculated kernels with psychophysiological findings reported in the recent literature. Then we will insert the connectivity kernels computed in Section \ref{sec:diffusion} in a neural population activity model \cite{Ermentrout1980}, regarding them as cortical facilitation patterns.

In section \ref{sec:numerical} we will use this population activity model equipped with the suitable connectivity kernels in two numerical simulations, comparing the results to recent phenomenological findings of the perception of contours in motion \cite{Rainville2005}, and to fMRI measurements of cortical neural activity related to motion perception \cite{Wu2011}.

\section{The geometry of spatio-temporal dynamics}\label{sec:geometry}

In this section we extend an approach introduced in \cite{CS} that amounts to model each V1 simple cell in terms of its receptive profile,  to interpret its action as a Gabor filtering, and to introduce a geometry of the space compatible with the properties of the output. In this paper we will consider each cell as sensitive to a local orientation and apparent velocity, that is the velocity orthogonal to a moving stimulus. The collected data lies on a five dimensional manifold $\M$ of space, time, orientation and velocity, a single cell being represented as a point on $\M$. The geometric structure of this manifold will be described in terms of a contact structure, that provides a constraint on the tangent space and on admissible connectivity among cells.

\subsection{Spatio-temporal receptive profiles}\label{sec:filtering}

It is known that the visual cortex decomposes the visual stimulus by measuring its local features. Local orientation and direction of movement have been the first visual features of neurons in V1 that have been studied. Receptive profiles (RPs) are descriptors of the linear filtering behavior of a cell and they can be reconstructed by processing electrophysiological recordings \cite{Ringach2002}. It has been shown that the spatial characteristics of these RPs can be modeled by 2-dimensional Gabor functions \cite{Daugman,Jones1987b}.

	\begin{figure*}
	\centering
	\includegraphics[width=\textwidth]{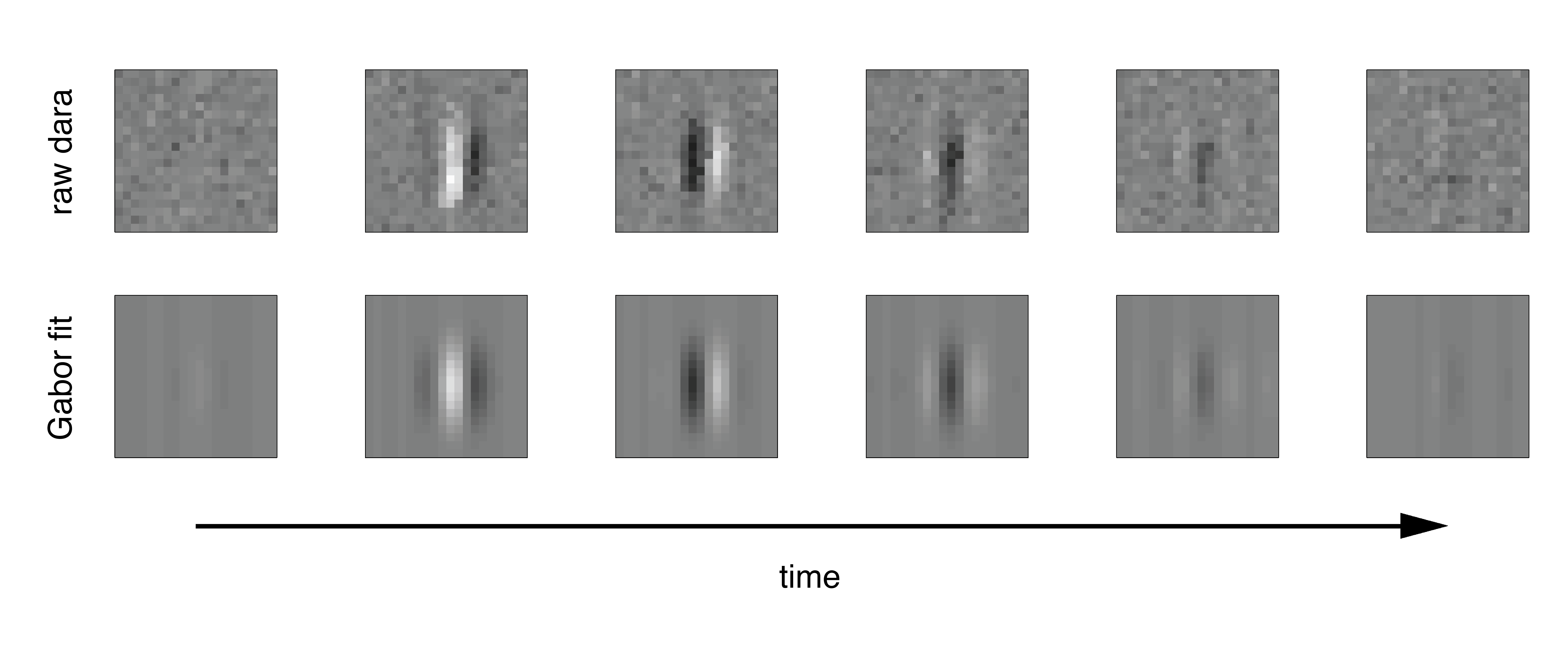}
	\caption{A spatio-temporal receptive profile taken from the data set studied in \cite{DeAngelis1993a} and \cite{CBS} (cell n.15 in the latter reference). In the top row, the sequence relative to the raw data, reconstructed using reverse correlation from electrophysiological recordings. In the bottom row, the 3-dimensional Gabor fit, visualized as 2D slices.}\label{fig:gaborseq}
	\end{figure*}
	
However a  large class of cells shows a very specific space-time behavior in which the spatial phase of the RP changes gradually as a function of time \cite{DeAngelis1993a}. Although many models have been proposed to reproduce these dynamics \cite{DeAngelis1993a,Adelson1985}, in \cite{CBS} it has been shown that a 3-dimensional sum-of-inseparable Gabor model can fit very well experimental data of both separable and inseparable RPs (Fig. \ref{fig:gaborseq}). Following this approach, we choose Gaussian Gabor filters centered at position $q = (q_1, q_2) \in \R^2$ on the image plane, activated around time $s \in \R$, with spatial frequency $p = (p_1,p_2) \in \R^2$, temporal frequency $\nu \in \R$, spatial width $\sigma_x \in \R^+$ (circular gaussians) and temporal width $\sigma_t \in \R^+$
\begin{equation}\label{eq:insep}
\RP(x,t) = e^{i(p\cdot (x - q) - \nu (t - s))} e^{-\frac{|x-q|^2}{\sigma_x^2}-\frac{(t-s)^2}{\sigma_t^2}}
\end{equation}
where we have used the abbreviations $z = (q_1, q_2, s) + i (p_1,p_2,\nu) \in \C^3$ and $\sigma = (\sigma_x, \sigma_t) \in \R^+ \times \R^+$. Also note that the functions (\ref{eq:insep}) correspond to the propagation  of a two dimensional plane wave within the activating window with phase velocity
\begin{equation}\label{eq:phasevel}
v = \frac{\nu}{|p|}.
\end{equation}
The variable $z = (q_1, q_2, s) + i (p_1,p_2,\nu) $ is canonically associated \cite{Folland} to the  phase space
\begin{equation}\label{eq:Sigma}
\R^6 = \{(q_1, q_2, s, p_1,p_2,\nu)\}
\end{equation}
endowed with the symplectic structure compatible with the complex structure of $\C^3$ of variable $z$, that is
\begin{displaymath}
\Omega = d\lambda = dp_1 \wedge dq_1 + dp_2 \wedge dq_2 - d\nu \wedge ds
\end{displaymath}
where we have denoted with $\lambda$ the Liouville form
\begin{equation}\label{eq:liouvilleform}
\lambda = p_1 dq_1 + p_2 dq_2 - \nu ds.
\end{equation}

\begin{figure*}
\centering
\includegraphics[width=\textwidth,height=.35\textwidth]{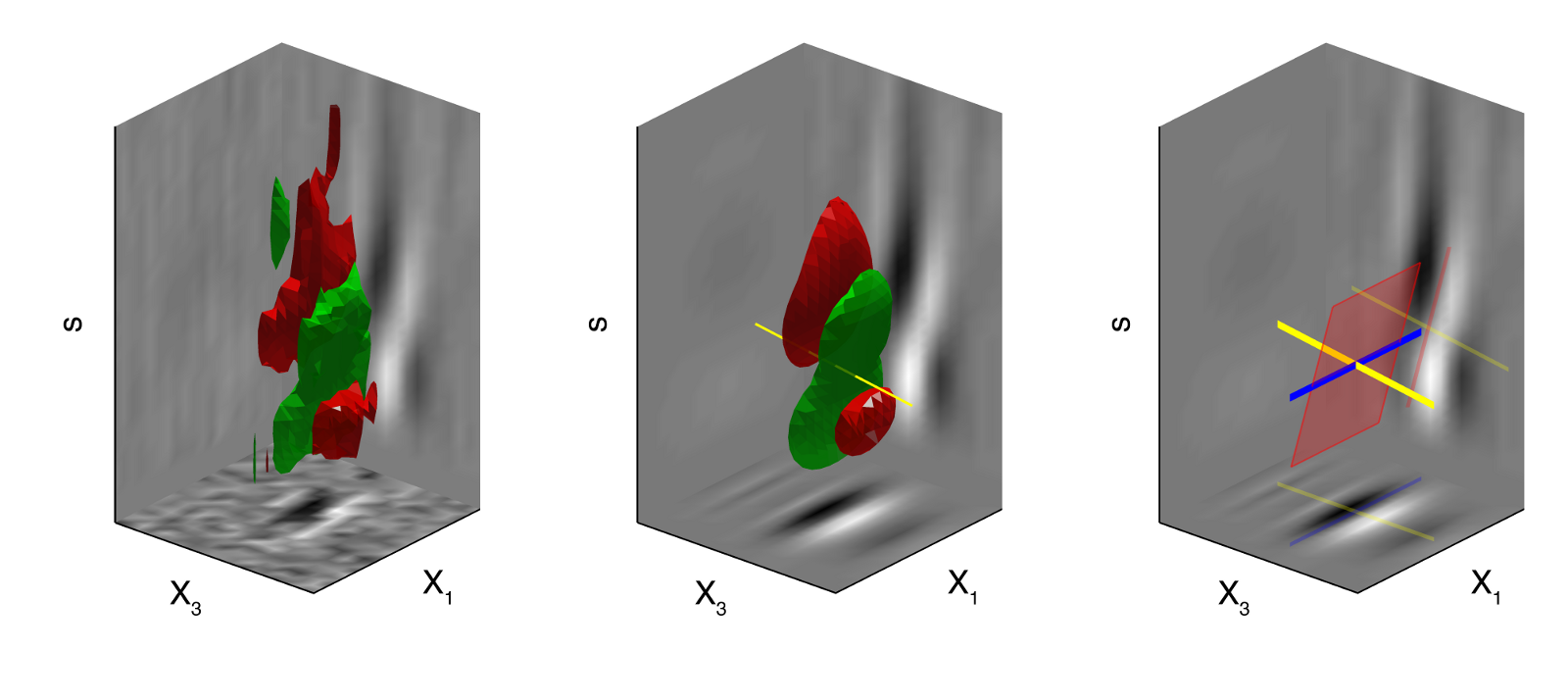}
\caption{Visualization of the one-form $\omega$ with respect to a cell RP. In the  left and center subplots are shown the isosurfaces of the RP already seen in Fig. \ref{fig:gaborseq}. In the left plot we visualize the reconstructed data, while in the middle plot we visualize its Gabor fit. The orientation of the yellow line within the fitted data represents the orientation of the wavefront of the Gabor filter and coincides with the one-form $\omega = \cos\theta dq_1+\sin\theta dq_2 - v ds$, or the dual vector $ \vec{X}$. In the right plot the semi-opaque red plane represents the horizontal tangent plane in the $R^3$ spatio-temporal space, generated by $ \vec{X}_5$ and $ \vec{X}_1$ (blue line).}  \label{fig:gaborform}
\end{figure*}

Since we are mainly interested in the geometry of level lines, our analysis can be restricted to the information captured by filters possessing a fixed central spatial frequency $|p|$. This amounts to fixing a scale of oscillations, hence disregarding the harmonic content of the filtering by focusing on the features of orientations and velocity. We then obtain the reduced 1-form
\begin{equation} \label{eq:contactform}
\omega = \frac{\lambda }{ |p|} = \cos\theta dq_1 + \sin\theta dq_2 - v ds 
\end{equation}
where $p= |p|(\cos\theta, \sin\theta)$. This form is defined on the spatio-temporal phase space with fixed frequency, that is the 5-dimensional manifold
\begin{displaymath}
\M = \R^2 \times \R^+ \times S^1  \times \R^+ = \{\eta=(q_1,q_2,s,\theta,v)\} ,
\end{displaymath}
and it is associated to every Gabor filter as shown in Fig. \ref{fig:gaborform}. At this level, time $s$ is introduced as a base variable with the same role of $\q$. The dual variable of $s$ is $v$, which has the same role of $\theta$, both being the engrafted variables with respect to time and space. To every point of $\M$ corresponds univocally a Gabor filter whose parameters are the coordinates of the point itself.

\subsection{Admissible tangent space as constraint on the connectivity on $\M$}\label{sec:contact}

We will model the connectivity between points in the space $\M$ in terms of admissible tangent directions of $\M$ itself. From the geometric point of view, the presence of the 1-form (\ref{eq:contactform}) is equivalent to the choice of a vector field with the same coefficients as $\omega$ with respect to the canonical basis $\{\frac{\partial}{\partial q_1}, \frac{\partial}{\partial q_2}, \frac{\partial}{\partial s}, \frac{\partial}{\partial \theta}, \frac{\partial}{\partial v}\}$, that is its dual vector (see Fig. \ref{fig:gaborform}): 
\begin{displaymath}
\vec{X} =(\cos\theta, \sin\theta,-v,0,0).
\end{displaymath}
The kernel of $\omega$, denoted $\ker \omega$ (or $\omega=0$) is the space of vectors orthogonal to $\vec{X}$. A basis of this space is constituted of the so-called horizontal or admissible vectors 
\begin{equation}\label{eq:vectorfields}
\begin{array}{l}
\vec{X}_1 = ( -\sin\theta,  \cos\theta, 0, 0, 0) \, , \ \vec{X}_2 = (0,0,0,1,0)\\
\vec{X}_4 = (0,0,0,0,1) \, , \ \vec{X}_5 = ( v \cos\theta, v \sin\theta, 1,0,0)
\end{array}
\end{equation}
so that
\begin{equation}\label{eq:hyperplanes}
\ker \omega = \vsp{\vec{X}_1, \vec{X}_2, \vec{X}_4, \vec{X}_5}
\end{equation}
that defines the {\it{horizontal tangent space}}. In particular we nota that the Euclidean metric on the horizontal planes makes the vector fields $X_i$ orthogonal.

% In the sequel we will assume that on the horizontal plane is defined a metric, which makes the vector fields $X_i$ orthogonal.

The only manifolds (curves or surfaces) admissible in this space are the ones whose tangent vectors are linear combinations of the horizontal vectors (\ref{eq:vectorfields}). The neural connectivity between the receptive fields at different points in the $\M$ space will be defined in terms of these vectors in Section \ref{sec:diffusion}. 

For reader convenience we compute here explicitly all the non-zero commutation relations between the vector fields \footnote{To every 
 vector field $\vec{X}= (a_1, a_2, a_3, a_4, a_5)$ we can associate a directional derivative 
$$X =a_1\partial_{q_1}+ a_2\partial_{q_1} + a_3 \partial_{s}+ a_4 \partial_{\theta}+ a_5 \partial_v$$ with the same coefficients. Then we will call commutator of $X$ and $Y$ $$[X,Y] = XY-YX.$$ We say that $X$ and $Y$ commute if $[X,Y]=0$. Note that partial derivatives always commute, while directional derivatives in general do not.  It is important to note that, even though the commutator is expressed formally as a second derivative, it is indeed a first derivative, so that it has an associated vector. Hence we can as well define the commutator between the vectors $\vec{X}$ and $\vec{Y}$ as the vector associated to $XY-YX$. }: 
\begin{eqnarray}
\ [\vec{X}_1 , \vec{X}_2] \ = & \vec{X}_3 &= (\cos\theta, \sin\theta, 0,0,0) \label{eq:commutglue} \\
\ [\vec{X}_2 , \vec{X}_3] \ = &\vec{X}_1& \nonumber \\
\ [\vec{X}_4 ,\vec{X}_5] \ = & \vec{X}_3& \label{eq:commutinside} \\
\ [\vec{X}_2 , \vec{X}_5] \ = &v\vec{X}_1. & \nonumber
\end{eqnarray}

\begin{figure*}
\centering
\subfloat[]{\includegraphics[width=.48\textwidth]{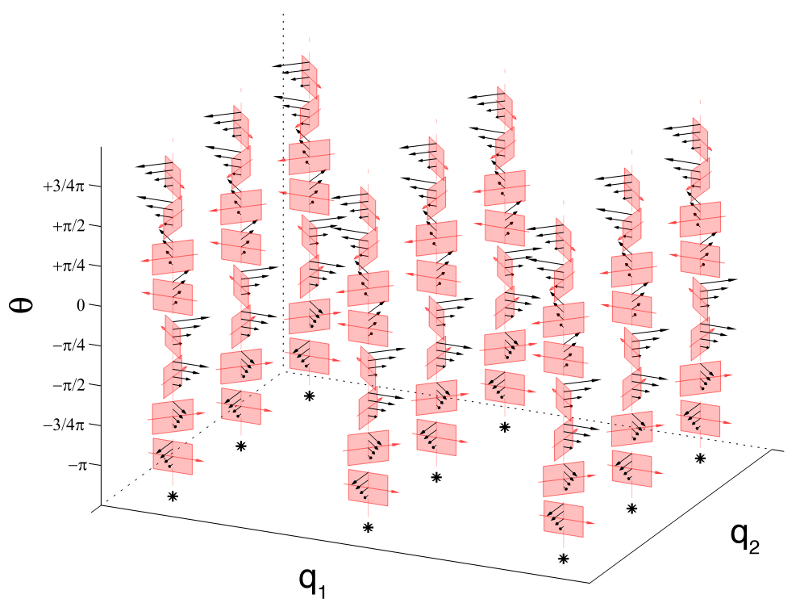}} \quad
\subfloat[]{\includegraphics[width=.48\textwidth]{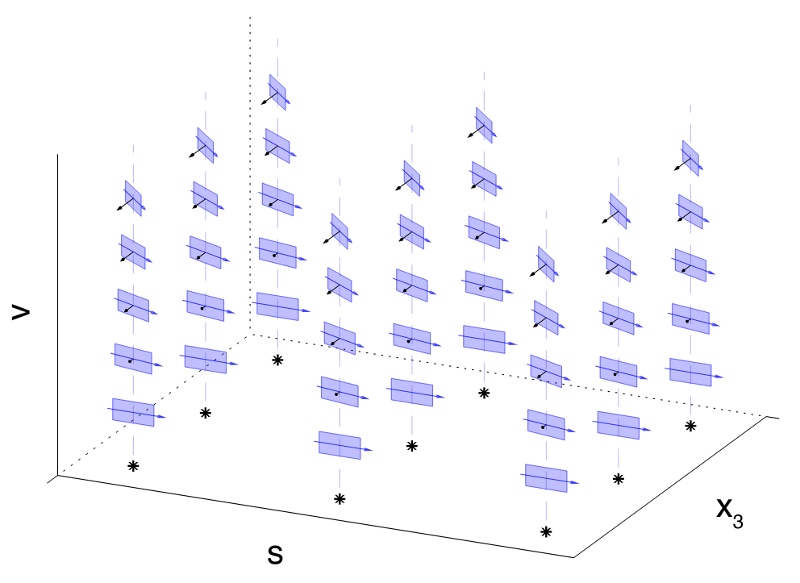}}
\caption{Visualization of the admissible tangent space of $\M$, generated by the vectors $\vec{X}_1, \vec{X}_2, \vec{X}_4,\vec{X}_5$. In a) we show the projection in the $(q,\theta)$ space. Red planes are the sections of the horizontal planes, generated by $\vec{X}_1$ and $\vec{X}_2$. Red arrows indicate the $\vec{X}_1$ direction, black arrows the $v\vec{X}_3$ one. In b) we depict the projection in the $(x_3, s, v)$ space. The blue planes are the sections of the  horizontal planes, spanned by $ \vec{X}_4$ and $ \vec{X}_5$  The black arrows represents the $v \vec{X}_3$ as before. }\label{fig:contactform}
\end{figure*}

In Fig. \ref{fig:contactform} we have depicted the structure of the tangent fields (\ref{eq:vectorfields}). For visualization purposes, we show the fields in two different representations.
In Fig. \ref{fig:contactform}a we visualize the structure restricted to the spatial and engrafted variables $(q,\theta,v)$, spanned by  $\vec{X}_1$ and $\vec{X}_2$. The tilting of the planes is due to the non commutative relation (\ref{eq:commutglue}). In Fig \ref{fig:contactform}b we show the spatio-temporal structure restricted to the variables $(x_3, s, v)$, where $x_3 = q_1 \cos \theta + q_2 \sin \theta$. Also in this figure the tilting of the planes spanned by $\vec{X}_4$ and $\vec{X}_5$ is due  to the non vanishing  commutator condition (\ref{eq:commutinside}). It is worth noting that in this setting, the propagation is permitted only along the horizontal planes, while it is forbidden in their orthogonal direction $\vec{X}$.
Let us note that 
$\vec{X}_3 = [\vec{X}_1 , \vec{X}_2]$
is linearly independent of the horizontal tangent space at every point. Hence $\{\vec{X}_1, \vec{X}_2, \vec{X}_4, \vec{X}_5\}$ together with their commutators span  the whole space at every point. 
This is the so called  H\"ormander condition, which will have a crucial role in studying the property of the space. The same condition can be expressed in terms of properties of the form $\omega$, which is called {\it contact form} \cite{Geiges} since  $\omega \wedge d\omega$ is never zero being the volume form of the space. 

We also note that we can define on $\M$ a smooth composition law
% \begin{equation}\label{eq:composition}
% \begin{array}{l}
% (\qv,s,\theta,v) \odot (\qp,s',\theta',v')\\
% = (R_\theta (\qp + \textstyle{\binom{v}{0}} s') + \qv, s' + s,\theta' + \theta, v' + v)
% \end{array}
% \end{equation}
\begin{equation}\label{eq:composition}
\begin{array}{l}
(q,s,\theta,v) \odot (q',s',\theta',v')\\
= (R_\theta (q' + \textstyle{\binom{v}{0}} s') + q, s' + s,\theta' + \theta, v' + v)
\end{array}
\end{equation}
where $R_\theta$ is a counterclockwise rotation of an angle $\theta$, that is such that the vector fields $\{\vec{X}_1, \vec{X}_2, \vec{X}_3, \vec{X}_4, \vec{X}_5\}$ are left invariant with respect to this law. This implies that also the admissible curves of the structrue and its kernels will be invariant. The manifold $\M$ together with this composition law can also be identified with a subset of the Galilei group (see Appendix).

\subsection{The output of the receptive profiles}\label{sec:output}
The output of the Gabor filters selects a set of points in $\M$ corresponding to specific features of the image. Since the Gabor filters are always connected in terms of the vectors (\ref{eq:vectorfields}), also their output 
will inherit this structure, and will be concentrated around an admissible surface. 

The energy output of a cell with receptive profile $\RP$ in presence of a spatio-temporal stimulus $f(x,t)$ is given by
\begin{equation}\label{eq:filtering}
\begin{array}{rcl}
F(q,s, \theta, v) & \doteq & \left| \langle \RP, f\rangle \right|^2 \\
& = & \Big| \displaystyle{\int_{\R^3} dx \, dt \, e^{-i |p|((x - q)_1 \cos\theta + (x - q)_2 \sin\theta -  v (t - s))}}\\
&& \hspace{48pt}e^{-\frac{|x-q|^2}{\sigma_x^2}-\frac{(t-s)^2}{\sigma_t^2}} f(x,t) \Big|^2.
\end{array}
% \begin{array}{l}
% F(q,s, \theta, v) \doteq \left| \langle \RP, f\rangle \right|^2 \\
% = \Big| \displaystyle{\int_{\R^3} dx \, dt \, e^{-i |p|((x - q)_1 \cos\theta + (x - q)_2 \sin\theta -  v (t - s))}}\\
% \hspace{53pt}e^{-\frac{|x-q|^2}{\sigma_x^2}-\frac{(t-s)^2}{\sigma_t^2}} f(x,t) \Big|^2.
% \end{array}
\end{equation}

We note that taking the square modulus in (\ref{eq:filtering}) disregards the phase of the corresponding linear filtering, which for many applications (see e.g. \cite{Duits2011}) is crucial. In this case however all the relevant information is encoded in this energy model. In particular the geometric quantities $\theta$ and $v$ are encoded in the points of maxima of the energy. To see why, we recall that the analogous of formula (\ref{eq:filtering}) with  purely spatial Gabor filters was studied in \cite{CS}, where the filtering output $F_0$ was a function of the variables $\q$ and $\theta$ alone (see also \cite{HP}). In order to outline the goemetric meaning of the lifting, we consider its action on level lines which are smooth. In that case, the output $F_0$ to a stimulus $f_0(x)$ takes its maximum around a value $\theta^* = \theta^*(\q)$
\begin{displaymath}
\max_{\theta} F_0 (\q,\theta) = F_0 (\q,\theta^*) 
\end{displaymath}
where $(\cos \theta^*, \sin \theta^*)$  identifies the orientation orthogonal to the level lines of $f_0(x)$. 
As a consequence, in \cite{CS} it is proved that the level lines of $f_0$ are lifted to curves admissible in the sense that their tangent vector lies in the kernel of the form 
\begin{displaymath}
\omega_3 = \cos\theta dq_1 + \sin\theta dq_2
\end{displaymath}
which is a contact form that can be obtained by restricting (\ref{eq:contactform}) to the space $(\q,\theta)$.

Since here we take into account time, then the output $F$ is a function defined on $\M$. In perfect analogy with the lower dimensional case, and denoting with $(\q, s(\q))$ a level set of the stimulus $f(x,t)$, the output $F$ takes its maximum around the values $\theta^* = \theta^*(\q,s(\q))$ and $v^* = v^*(\q,s(\q))$
\begin{displaymath}
\max_{\theta , v} F (\q,s(\q),\theta,v) = F (\q,s(\q),\theta^*,v^*)
\end{displaymath}
such that the vector $(\cos\theta^*, \sin\theta^*, v^*)$ is orthogonal to the level set of $f$. The vector $(\cos\theta^*,\sin\theta^*)$ is orthogonal to the spatial level line, and the scalar $v^*$ represents the apparent velocity in this direction.

\begin{figure}
\centering
\includegraphics[width=.48\textwidth]{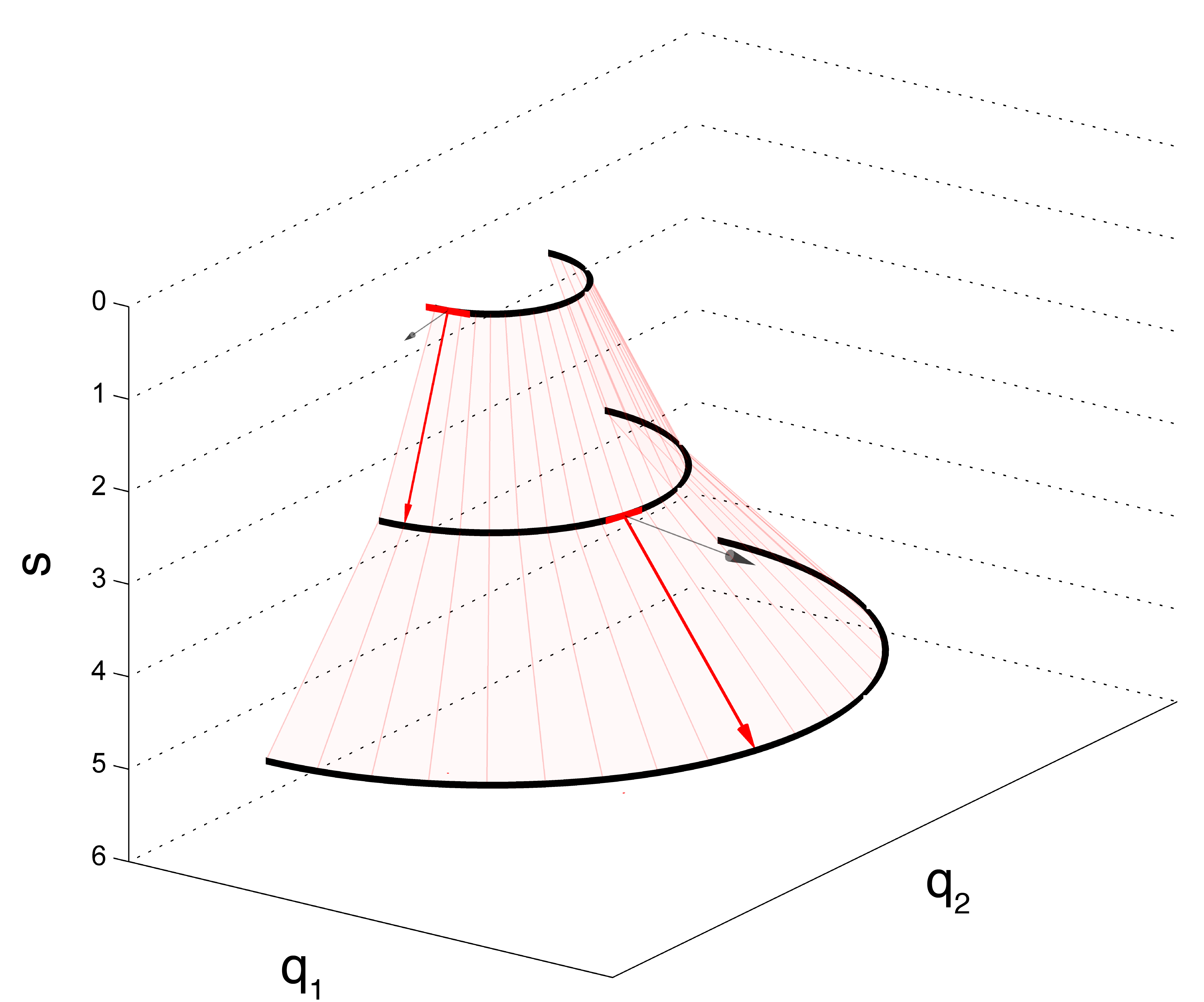} 
\caption{A schematical representation of the level set of a function $f$ (in this case, an expanding contour) lifted as an admissible surface $\Sigma$ on the 5D structure. Red and black arrows indicate, respectively, the local direction of the vector fields $\vec{X_5}$ and $\vec{X_3}$ over two points in space-time.}\label{fig:lslifting}
\end{figure}

The surface
\begin{displaymath}
\Sigma=\Big\{\big(\q, s(\q), \theta^*, v^*\big)\Big\}
\end{displaymath}
is the 5D lifting of the level set of $f$, and the orthogonality condition implies that it is admissible, in the sense that its tangent vectors lie in the horizontal space, kernel of the form $\omega$ (see Fig. \ref{fig:lslifting}). 

This shows that such level sets of $f$ define a contact structure on the manifold, so the study of such geometric features induces to endow $\M$ with the corresponding constraint on the tangent bundle that can be equivalently interpreted as a sub-Riemannian constraint \cite{Montgomery,CS,BDRS}.

\section{Curves and kernels of connectivity}\label{sec:diffusion}

The lifted points of the spatio-temporal stimulus are connected by admissible integral curves. However we will see that not all admissible curves can be considered lifted ones, and we will describe which ones have this property. These curves will represent the association fields in space-time, analogously to the association fields of Field Hayes and Hess\cite{Field1993}  in the pure spatial case, and are reminiscent of the classical Gestalt concept of \emph{good continuation}. We will discuss in Section \ref{sec:truelife} the compatibility of this structure with the cortical functional architecture. A probabilistic version of the connectivity field will be provided on the basis of Fokker Planck equations first introduced by D. Mumford in the spatial case \cite{Mumford}, and interpreted as model of connectivity in Lie groups in several works such as \cite{Petitot1999, WilliamsZweck, AugustZucker2003, Duits2010a, Duits2010b, SCS}. Here we will extend this stochastic approach to the space-time contact structure.

\subsection{Generators of lifted curves}\label{sec:legendrian}
We have already seen that the output of RP filtering  is concentrated around admissible submanifolds. Nevertheless not all admissible submanifolds can be lifting of images, since lifting are graphs of the functions 
$ \theta^*$ and $ v^*$. For example the plane generated by the vectors $\vec{X}_2$ and $\vec{X}_4$ cannot be recovered by lifting. Hence we will study the possible linear combinations of vector fields $\{\vec{X_i}\}$, which can be tangent to lifted level lines. An $\omega$-admissible curve $\gamma \subset \M$ is characterized by $\dot{\gamma} \in \ker\omega$, i.e.
\begin{displaymath}
\dot{\gamma} = a_1 \vec{X}_1(\gamma) + a_2 \vec{X}_2(\gamma) + a_4 \vec{X}_4(\gamma) + a_5 \vec{X}_5(\gamma)
\end{displaymath}
with $a_i$ not necessarily constant. Lifted curves  depend on the vectors $ \vec{X}_1$ or $ \vec{X}_5$, which are tangent to the base space $(q, s)$, and their linear combination. To simplify the problem, we will consider separately two special cases of particlar interest for the model: the limit cases of contour motion detected at a fixed time $(a_5=0)$, and motion of a point in time $(a_1=0)$. The first one is described by integral curves of the vector $\vec{X_1}$ and of the generators of the engrafted variables 
$\vec{X_2}$ and $\vec{X_4}$:
\begin{equation}\label{eq:contour}
\dot{\gamma}_{3} = \vec{X}_1(\gamma) + k \vec{X}_2(\gamma) + c \vec{X_4}(\gamma).
\end{equation}
These curves lie in the section of the contact structure depicted in Figure \ref{fig:contactform}a, where $k$ represents Euclidean curvature and $c$ is the rate of change of local velocity along the curve. On the other hand, the motion of a point in time can be described as 
\begin{equation}\label{eq:trajectory}
\dot{\gamma}_{T} = \vec{X}_5(\gamma) + w \vec{X}_2(\gamma) + a \vec{X}_4(\gamma)
\end{equation}
that are suitable to describe spatio-temporal trajectories of points. The coefficient $w$ of the direction $\vec{X}_2$ is the angular velocity, and the coefficient $a$ of $\vec{X}_4$ is the tangential acceleration.

\subsection{Curves and kernels for contours in motion}\label{sec:nt}

\subsubsection{Moving contours as deterministic integral curves}

We consider now the geometry  of moving contours at a fixed time. 
These are generated by the integral curves (\ref{eq:contour}) so to satisfy the system of equations
\begin{equation}\label{eq:integralcurvesV}
\left\{
\begin{array}{rcl}
\dot\gamma(t) & = & \vec{X}_1(\gamma) + k \vec{X}_2(\gamma) + c \vec{X}_4(\gamma)\vspace{4pt}\\
\gamma(0) & = & \xi_0
\end{array}
\right.
\end{equation}
with variable coefficients $k$, $c$. In particular, due to condition (\ref{eq:commutglue}) the set of points reacheable by piecewise constant integral curves of this type (for an explicit expression see e.g. \cite{CS}) starting from a fixed point $\xi_0 = (q_0,s_0,\theta_0,v_0) $ is the space 
\begin{equation}\label{enne}
\N = \{ (\q, s_0, \theta, v)\}
\end{equation}
where $s_0$ is fixed, so the space $\N$ can be identified with the space of points  $\xi = (\q,   \theta, v)$.

\begin{figure}
\centering
\includegraphics[width=.48\textwidth]{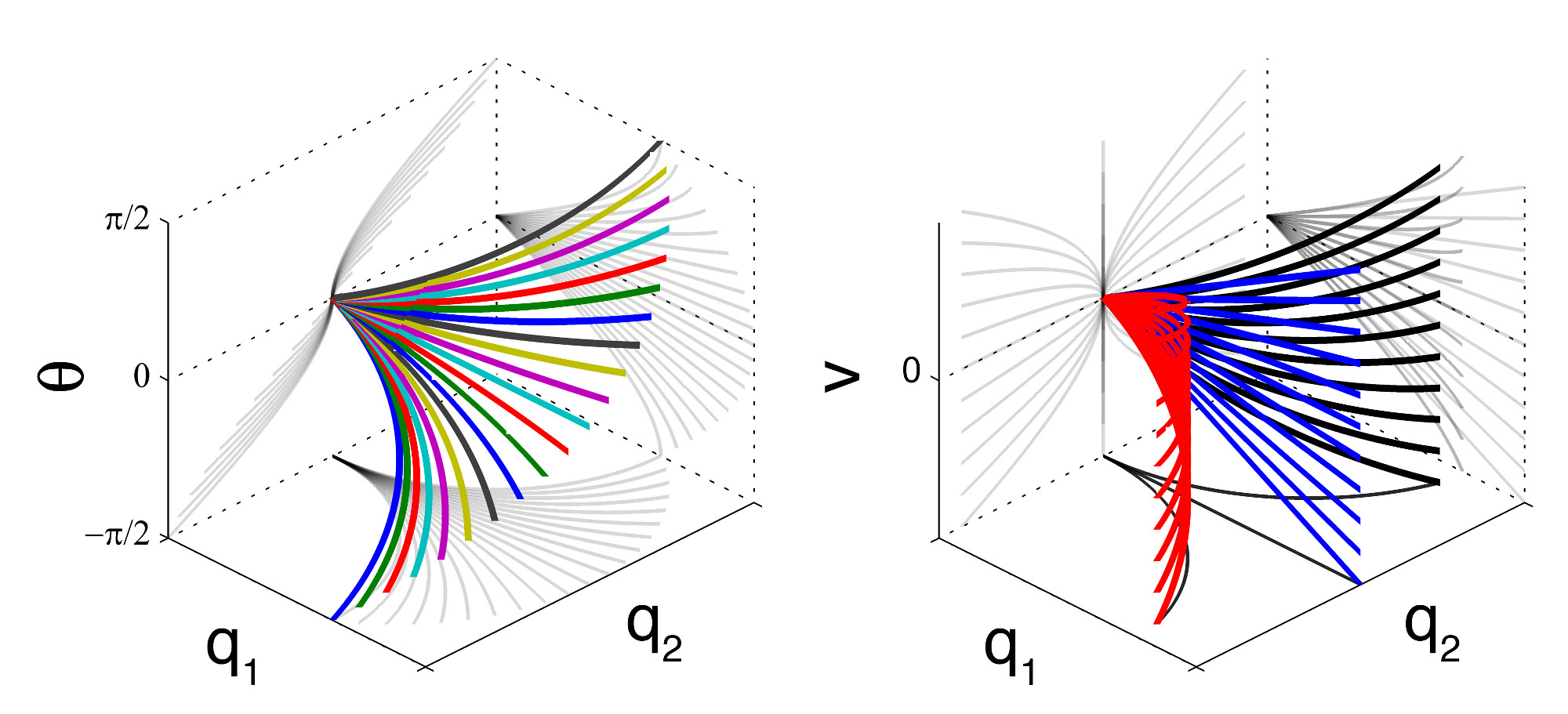}
\caption{The subset of horizontal curves with constant coefficients pertaining to $\Sigma_0(\xi_0)$ by varying parameters $k$ and $c$, with initial velocity set to $0$. Left plot: projections over the variables $(q_1, q_2, \theta)$ with different values for $k$, and same non-null value for $c$. Right plot: projections over the variables $(q_1, q_2, v)$ relative to some negative and positive values of $c$. The curves with $k = 0$ (blue ones in the right plot) are straight lines, as the direction of spatial propagation never changes.}\label{fig:horcurves1}
\end{figure}

We also denote with $\Sigma_0(\xi_0)$ the set of points reached by the fan of curves solution to the system (\ref{eq:integralcurvesV}) with constant coefficients $(k,c)$, depicted in Fig. \ref{fig:horcurves1}.

As we did with (\ref{eq:composition}) we can define on $\N$ a smooth composition law
% \begin{equation}\label{eq:composition}
% \begin{array}{l}
% (\qv,s,\theta,v) \odot (\qp,s',\theta',v')\\
% = (R_\theta (\qp + \textstyle{\binom{v}{0}} s') + \qv, s' + s,\theta' + \theta, v' + v)
% \end{array}
% \end{equation}
\begin{equation}\label{eq:compositionV}
(q,\theta,v) \boxdot (q',\theta',v') = (R_\theta q' + q, \theta' + \theta, v' + v)
\end{equation}
where $R_\theta$ is a counterclockwise rotation of an angle $\theta$, that is such that the vector fields $\{\vec{X}_1, \vec{X}_2, \vec{X}_3, \vec{X}_4\}$ are left invariant with respect to this law. Since this is a group law, the manifold $\N$ together with this (\ref{eq:compositionV}) provides the Lie group $\R^2 \ltimes S^1 \times \R$, that is the direct product of the $SE(2)$ group and the group of the reals with addition.

\subsubsection{Stochastic kernel}\label{sec:stochasticX2X4X1}

\begin{figure}
\centering
\subfloat[]{\includegraphics[width=.48\textwidth]{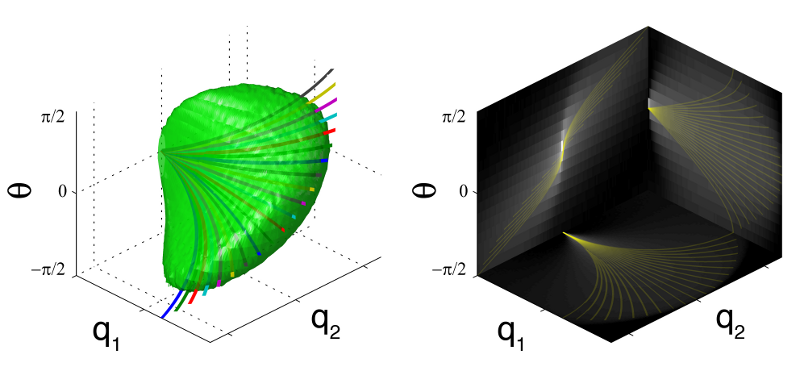}\label{fig:kernels1a}} \quad
\subfloat[]{\includegraphics[width=.48\textwidth]{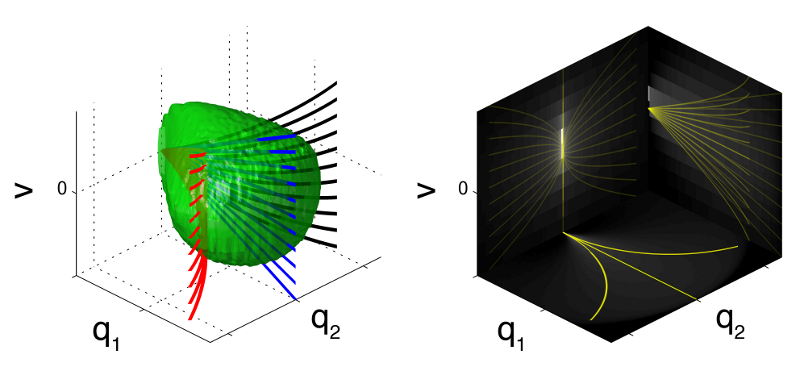}\label{fig:kernels1b}}
\caption{Horizontal curves and stochastic kernels for the integration of contours in motion. Left: the isosurface plot of the kernel $\Gv$ (isovalue: $0.002 \max(\Gv)$) is superimposed over the subsets of curves (\ref{eq:integralcurvesV}) with constant coefficients already shown in Fig. \ref{fig:horcurves1}. Right: the kernel projections of $\Gv$ relative to the variables a) $(q_1, q_2, \theta)$ and b) $(q_1,q_2,v)$ are plotted in gray under the projections of the curves (in yellow). This kernel was obtained by calculating the evolution of $1 0^6$ stochastic paths.}\label{fig:kernels1}
\end{figure}

Let us consider a probabilistic counter part of trajectories (\ref{eq:integralcurvesV}). The vector field $X_1$ expresses a derivative in the direction of a variable coded on the retina, while the vectors $X_2$, $X_4$ express derivatives in the direction of an eingrafted variable. Due to their different role, we will consider the following vector-valued stochastic process
\begin{equation}\label{eq:brown1}
\left\{
\begin{array}{rcl}
d\gamma & = & \vec{X}_1(\gamma) dt + \vec{X}_2(\gamma) dW_1 +\vec{X}_4(\gamma) dW_2\vspace{0.4pt}\\
\gamma(0) & = & \xi_0
\end{array}
\right.
\end{equation}
where $W = (W_1, W_2)$ is a two dimensional Brownian motion. The distribution of these stochastic curves is mostly concentrated around the surface $\Sigma_0(\xi_0)$. The density of points reached by this stochastic kernel is then a candidate to implement the mechanism of association fields (see Fig. \ref{fig:kernels1}). This approach generalizes the approach of random paths introduced in \cite{Mumford, Williams} for the $SE(2)$ problem.

If we call $\rv(\xi,t|\xi_0,0)$ the density of points of $\N$ reached at the value $t$ of the evolution parameter by the sample paths of the process (\ref{eq:brown1}), then $\rv$ can be obtained as the fundamental solution
\begin{displaymath}
(\partial_t + \Lv)\rv(\xi,t|\xi_0,0) = \delta(\xi-\xi_0)\delta(t)
\end{displaymath}
where $\Lv$ is the Fokker-Planck operator 
\begin{equation}\label{eq:Lgv}
\Lv \doteq X_1 - \kappa^2 X_2^2 - \alpha^2 X_4^2 
\end{equation}
containing a diffusion over the fiber variables $\theta$ and $v$ and a drift in the base variables, and with $X_i$ and $X_i^2$ we denote first and second order derivatives in direction $\vec{X_i}$. Since $\Lv$ contains a set of vector fields that generates the whole tangent space of $\N$, by H\"ormander theorem this operator is hypoelliptic \cite{Hormander}. This means that $\rv$ is non-null in all $\N$ for any $t>0$, even if the operator $\Lv$ contains only 3 linearly independent fields. Moreover, since it defines the Fokker-Planck equation for the stochastic process (\ref{eq:brown1}), $\rv$ is indeed a (conditional) probability density on $\N$, evolving with the parameter $t$.

In order to characterize each point of $\N$ in terms of the density of paths (\ref{eq:brown1}) that reach it, independently of the value of the evolution parameter, we need a notion corresponding to the fan $\Sigma_0(\xi_0)$, expressed in terms of $\rv$. The density of points reached \emph{at any value of the evolution parameter} by the stochastic dynamics (\ref{eq:brown1}) is given by
\begin{equation}\label{eq:densityg0}
\Gv(\xi|\xi_0) \doteq \int_{\R} \rv(\xi,t|\xi_0,0) dt\ . 
\end{equation}
This derived quantity is actually the fundamental solution of the operator $\Lv$, so explicitly we have
\begin{equation}\label{eq:FPgv}
\Big(X_1 - \left(\kappa^2 \partial^2_\theta + \alpha^2\partial^2_v \right) \Big) \Gv(\xi|\xi_0) = \delta(\xi-\xi_0)
\end{equation}
and since the vector fields involved in equation (\ref{eq:Lgv}) are left invariant with respect to the group law (\ref{eq:compositionV}), the solution (\ref{eq:densityg0}) possesses the symmetry
\begin{equation}\label{eq:symmetryV}
\Gv(\xi|\xi_0) = \Gv(\xi_0^{-1}\boxdot \xi|0)\quad \forall \ \xi,\xi_0 \in \N\ .
\end{equation}
% Analogously, the fan $\Sigma(\eta_0)$ of integral curves (\ref{eq:associationT}) with constant coefficients starting from a point $\eta_0$ can be obtained from the one starting from $0$ by
% \begin{displaymath}
% \Sigma(\eta_0) = \eta_0^{-1_L}\odot \Sigma(0)\ .
% \end{displaymath}

We recall that in \cite{DvA} several representations of the exact solution to a problem analogous to (\ref{eq:FPgv}) were presented. In this work we will deal with the numerical implementation of the fundamental solution $\Gv(\xi|\xi_0)$, developed with standard  Markov Chain Monte Carlo methods (MCMC) \cite{MCMC}. This is done by generating random paths obtained from numerical solutions of the system (\ref{eq:brown1}) and averaging their passages over discrete volume elements, and appears suitable to treat also the more involved case of subsequent equation (\ref{eq:FPg0}).
An example is shown in Fig. \ref{fig:kernels1}, where an isosurface plot of the kernel is plotted over the integral curves (\ref{eq:integralcurvesV}) with constant coefficients, already depicted in Fig. \ref{fig:horcurves1}. By comparing such numerical approximations with \cite{DvA}, we can confirm the accuracy of this method. Moreover, from the figure it can be seen that the probability density is concentrated around the surface $\Sigma_0(\xi_0)$ and decays rapidly away from it.  This is reasonable since for this kind of hypoelliptic operators one can generally obtain estimates for the fundamental solution in terms of exponential decay with respect to a geodesic distance computed with respect to the minimal set of vector fields that, together with their commutators, span the whole tangent space \cite{RS}. Such a distance is anisotropic, and its balls are squeezed in the directions of the commutators \cite{NSW}. Here the vector fields involved are $\{X_1, X_2, X_4\}$, hence the concentration of the fundamental solution around the set $\Sigma_0$ defined by their integral curves follows by the sharper decay in the direction of the commutators. We also note that in this particular case, due to the availability of exact solutions and estimates from \cite{DvA}, this behaviour can also be checked directly.

\subsection{Curves and kernels for point trajectories}\label{sec:sR}

\subsubsection{Point trajectories as deterministic integral curves}

The trajectory-type curves introduced in (\ref{eq:trajectory}) are solutions $\gamma: \R \to \M$ to the system of ordinary differential equations
\begin{equation}\label{eq:associationT}
\left\{
\begin{array}{rcl}
\dot\gamma(t) & = & \vec{X}_5(\gamma) + w \vec{X}_2(\gamma) + a \vec{X}_4(\gamma)\vspace{4pt}\\
\gamma(0) & = & \eta_0
\end{array}
\right.
\end{equation}
for given initial point $\eta_0 = ({q_1}_0,{q_2}_0,s_0,\theta_0,v_0)$. In general the coefficients $w, a$ need not to be constant, but we can have a local approximation of any curve in a neighborhood of the starting point if we consider the case of constant coefficients. In this model case when $w$ and $a$ are not zero, (\ref{eq:associationT}) is explicitly solved by
\begin{displaymath}
\left\{
\begin{array}{rcl}
q_1(t) & = & {q_1}_0 + \varrho (\cos(\theta(t)-\phi) -  \cos(\theta_0-\phi)) +  \displaystyle{\frac{a}{w}} t \sin\theta(t) \vspace{4pt}\\
q_2(t) & = & {q_2}_0 + \varrho (\sin(\theta(t)-\phi) -  \sin(\theta_0-\phi)) -  \displaystyle{\frac{a}{w}} t \cos\theta(t)\vspace{4pt}\\
s(t) & = & s_0 + t\vspace{4pt}\\
\theta(t) & = & \theta_0 + w t\vspace{4pt}\\
v(t) & = & v_0 + a t
\end{array}
\right.
\end{displaymath}
where $\varrho = \frac{\sqrt{a^2 + v_0^2 w^2}}{w^2}$ and $\phi = \arctan\frac{v_0 w}{a}$, and a reasonable choice is to set $s_0 = 0$, in order to synchronize the evolution parameter $t$ with the time parameter $s$.

\begin{figure}
\centering
\includegraphics[width=.48\textwidth]{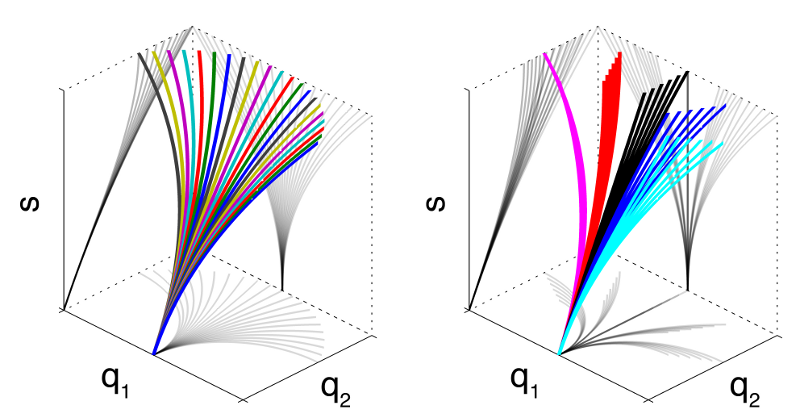}
\caption{Projections over the variables $(q_1, q_2, s)$ of the fan  of curves with constant coefficients, for different values of parameters $w$ and $a$ and a non-null initial velocity. Different values of $w$ are associated to different curve colors, while gray is used for projections over two dimensional planes. Left plot: curves with different values of $w$, and a fixed non-null value for $a$. Right plot: curves with different values of $a$.}\label{fig:horcurves2}
\end{figure}

The fan $\Sigma(\eta_0)$ of such curves is depicted in Fig. \ref{fig:horcurves2}. Each of them describes a motion on the plane $\q$ that, for small times, consists approximately of arcs of circles with radius $\varrho$, while for sufficiently large times corresponds to enlarging spirals around a slightly moving center. Their instantaneous acceleration is given by $\displaystyle{\binom{\ddot{q}_1}{\ddot{q}_2} = R_\theta \binom{a}{v w}}$, so that $a$ corresponds to
the tangential acceleration\footnote{By direct computation, $a$ is the time derivative of the modulus of the velocity
\begin{displaymath}
a = \frac{d}{dt} \sqrt{\dot{q}_1(t)^2 + \dot{q}_2(t)^2}\ .
\end{displaymath}} along direction $X_3$:
\begin{displaymath}
a = \ddot{q}_1 \cos\theta + \ddot{q}_2 \sin\theta
\end{displaymath}
while their  curvature is
\begin{displaymath}
\frac{\dot{q}_1\ddot{q}_2 - \dot{q}_2\ddot{q}_1}{(\dot{q}_1{}^2 + \dot{q}_2{}^2)^{\frac32}} = \frac{w}{v}\ .
\end{displaymath}
In particular, when $w = 0$ we obtain straight lines along the $\theta_0$ direction, while for $a = 0$ we obtain circular trajectories of radius $\displaystyle{\frac{v_0}{w}}$.

\subsubsection{Stochastic kernel}\label{sec:stochasticX2X4X5}

Similarly to what we have done in section \ref{sec:stochasticX2X4X1}, let us consider the vector-valued stochastic process
\begin{equation}\label{eq:brown2}
\left\{
\begin{array}{rcl}
d\gamma & = & X_5(\gamma) dt + X_2(\gamma) dW_1 + X_4(\gamma) dW_2\vspace{0.4pt}\\
\gamma(0) & = & \eta_0
\end{array}
\right.
\end{equation}
where $W = (W_1, W_2)$ is a two dimensional brownian motion. 

The density $\ro(\eta,t|\eta_0,0)$ of points of $\M$ reached at the value $t$ of the evolution parameter by the sample paths of the process (\ref{eq:brown2}), is the fundamental solution of the equation
\begin{equation}\label{eq:FPg0NH}
(\partial_t + \Lo)\ro(\eta,t|\eta_0,0) = \delta(\eta-\eta_0)\delta(t)
\end{equation}
where $\Lo$ is the Fokker-Planck operator 
\begin{equation}\label{eq:Lg0}
\Lo \doteq - \alpha^2 X_4^2 - \kappa^2 X_2^2 + X_5.
\end{equation}
Equation (\ref{eq:FPg0NH}) is not hypoelliptic, indeed its fundamental solution $\rho$ is concentrated on a submanifold of codimension 1 defined by the equation $t = s$, from the system (\ref{eq:brown2}). However it is still a Fokker-Planck equation, hence $\rho$ is nonnegative, and integrating the density $\ro$ with respect to the evolution parameter $t$ we derive the fundamental solution of the operator $\Lo$. Explicitly we have
\begin{equation}\label{eq:FPg0}
\Big(\partial_s - \left(\kappa^2 \partial^2_\theta + \alpha^2\partial^2_v - v X_3\right) \Big) \Go(\eta|\eta_0) = \delta(\eta-\eta_0)
\end{equation}
and we note in particular that this equation is now a hypoelliptic equation that is also the Fokker-Planck equation of a stochastic process defining a propagation in the physical time $s$.
Moreover, since the vector fields that consitute the operator (\ref{eq:Lg0}) are left invariant with respect to the composition law (\ref{eq:composition}), the solution (\ref{eq:densityg0}) possesses the symmetry
\begin{equation}\label{eq:symmetry}
\Go(\eta|\eta_0) = \Go(\eta_0^{-1_L}\odot \eta|0)\quad \forall \ \eta,\eta_0 \in \M
\end{equation}
where $\eta_0^{-1_L}$ stands for the left inverse with respect to (\ref{eq:composition}), see also Appendix.

\begin{figure}
\centering
\subfloat[]{\includegraphics[width=.48\textwidth]{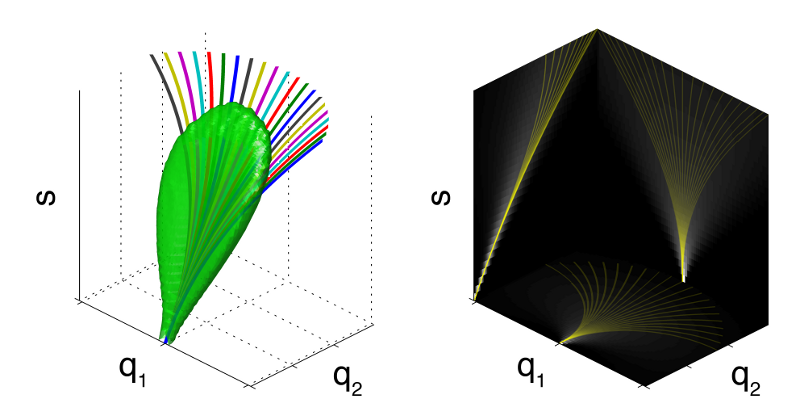}} \quad
\subfloat[]{\includegraphics[width=.48\textwidth]{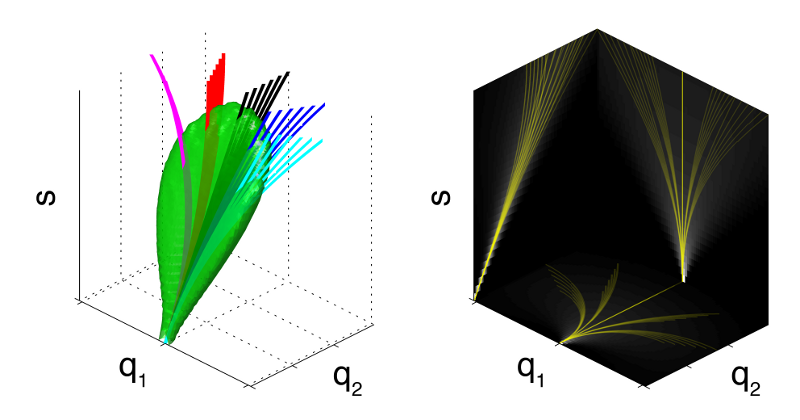}}
\caption{Horizontal curves and stochastic kernels for the integration of motion, with a non-null initial velocity value. As before, different values of $k$ are associated to different curve colors. Left: the isosurface plot of the kernel $\Go$ (isovalue: $0.001 \max(\Go)$) is superimposed over the curves (\ref{eq:associationT}) with constant coefficients already shown in Fig. \ref{fig:horcurves2}. Right: the projections of $\Go$ relative to the variables $(q_1, q_2, s)$ are plotted under the projections of the curves (in yellow). This kernel was obtained by calculating the evolution of $10^6$ stochastic paths.}\label{fig:kernels2}
\end{figure}

Similarly to what we have done in the previous section, we can obtain $\Go(\eta|\eta_0)$ at each point by solving numerically the system of stochastic differential equations (\ref{eq:brown2}) and applying standard Markov Chain Monte Carlo methods. An isosurface plot of the kernel is shown in Fig. \ref{fig:kernels2}, over the same integral curves of Fig. \ref{fig:horcurves2} and again we see the concentration around the fan $\Sigma(\eta_0)$.

The fundamental solution (\ref{eq:FPg0}) will be concretely applied in Section \ref{sec:numerical} as a facilitation field, toward the task of motion integration of trajectories.

\section{Neural propagation of boundaries and trajectories}\label{sec:truelife}

The curves (\ref{eq:integralcurvesV}) and (\ref{eq:associationT}) can be related to well-defined perceptual mechanisms. More precisely, we propose to consider them as association fields, in the sense of \cite{Field1993}, devoted to the integration of contour in motion and trajectories. The perceptual tasks of subjective boundary completion and motion integration have been widely studied by both psychologists and physiologists, and their underlying physiological explanation continues to be an open issue of discussion. The perceptual bias towards collinear stimuli has classically been associated to the long-range horizontal connections linking cells in V1 sharing similar preferences in stimulus orientation. This specialized form of intra-striate connectivity pattern is found across many species, including cats, tree shrews and macaques \cite{Kisvarday1997,Chisum2003}, the main difference being the specificity and the spatial extent of the connections. Furthermore, axons seem to follow the retinotopic cortical map anisotropically, with the axis of anisotropy being related to the orientation tuning of the originating cell \cite{Bosking1997}. These connections have already been modeled in \cite{CS} by means of a sub-Riemannian diffusion process over the orientation space $\R^2 \times S^1$.

Similarly to what happens for the integration of merely spatial visual information, the brain is also capable to easily predict stimulus trajectories \cite{Verghese2002}, and to group together boundary elements sharing similar motion or apparent motion paths \cite{Ledgeway2005,Rainville2005}. One possible explanation for these effects could be the existence of specialized facilitatory networks linking cells anisotropically and coherently with their preference in velocity and axis of motion direction. Supporting these assumptions, it has been found that the neural preferences in direction of movement are also structurally mapped in the cortical surface, with nearby neurons being tuned for similar motion direction \cite{Weliky1996}, and it has been shown that excitatory horizontal connections in the V1 of the ferret are strictly iso-direction tuned \cite{Roerig1999}. Furthermore, it is known that also extra-striate area MT/V5 is retinotopically organized, its horizontal connectivity pattern being highly structured, with connections reaching columns of cells tuned for similar orientation and direction preference anisotropically and asymmetrically \cite{Malach1997, Ahmed2012}. Moreover, striate and extra-striate cortical areas seem to cooperate, and surround modulation in V1 can be given by the connectivity patterns implemented in both areas by means of fast feedforward and feedback inter-areal projections \cite{Angelucci2002}.

Regarding the physiological correlates of motion integration, it was shown that the motion of an occluded object trajectory is significantly represented in the human brain by the same visual areas that process real motion \cite{Olson2004}. Moreover, a recent study based on electrophysiological recordings of the V1 of tree shrews showed some non-linear neural behaviors that are coherent with the phenomenological dynamics of motion integration \cite{Wu2011}. It is indeed a general assumption that some cortical area in the visual cortex is responsible for predicting future motion, a possible implementation being a specialized connectivity for spatio-temporal trajectory facilitation, that is different from the one in V1 responsible for contour integration \cite{Grzywacz1995, Verghese2000, Watamaniuk2005, Whitaker2008}.

Here, we do not speculate on the exact physiological origin of these psychophysiological findings. However, we want to show that the connectivity kernels $\Gv$ and $\Go$ arising from the geometry defined in Section \ref{sec:geometry} is capable to reproduce qualitatively some of the effects reported on the works that have been previously cited. We will embed this connectivity in the neural population activity model described in the next paragraph. Then, in Section \ref{sec:numerical} we will simulate the response of cortical visual neurons to artificial stimuli, comparing the results with some psychophysiological findings reported in the recent literature.

\subsection{Modeling neural activity}

The state of a population of cells can be characterized by a real-valued activity variable, which depends on the interaction of the feedforward input at different points, due to the cortical connectivity. The first population activity models are due to Wilson and Cowan \cite{Wilson1972}, and Ermentrout and Cowan \cite{Ermentrout1980}: 
\begin{equation}\label{EC}
\frac{\textrm{d}a(\eta,t)}{\textrm{d}t} = -a(\eta,t) + S \Big( c_f \int_\M \Gamma(\eta |\zeta)  a(\zeta,t) d\zeta  + F(\eta)\Big)
\end{equation}
where $a$ is the neural activity of the population, $F$ is the feedforward input (\ref{eq:filtering}), $\Gamma$ is the facilitation kernel, $c_f$ is the facilitation strength and $S$ is the sigmoidal function 
\begin{equation}\label{eq:v1thresholding}
S( \tau) = \frac{1}{1 + e^{ - \mu ( \tau - \beta) }}\ .
\end{equation}
In the stationary case a first order approximation of the solution of (\ref{EC}) is 
\begin{equation}\label{giac}
F_0(\eta) = S \Big( c_f \int_\M \Gamma(\eta |\zeta) S(  F(\zeta) ) d\zeta  + F(\eta)\Big)
\end{equation}
that is the activity formula that we will use in the experiments of Section \ref{sec:numerical}. Let's explicitly note that the term
\begin{equation}\label{eq:ec}
F_T(\eta)= S(F(\eta))
\end{equation}
represents the mean neural extra-cellular activity in response to a stimulus. 

The geometry of the functional architectures is contained in the kernel $\Gamma$, taking into account the deep structures of the connectivity space.  Then the overall probability of activation can be obtained by convolution of 
the activity with the kernel $\Gamma$, and the term
\begin{equation}\label{eq:cc}
P(\eta)=\int_\M \Gamma(\eta |\zeta) F_T(\zeta) d\zeta
\end{equation}
is the facilitation pattern resulting from the contribution of horizontal cortico-cortical or feedback inter-areal connectivities. 

We note that, since the Lebesgue measure $d\zeta$ on $\M$ is left invariant with respect to the composition law (\ref{eq:composition}), and due to the symmetry (\ref{eq:symmetry}), then (\ref{eq:cc}) has the structure of a convolution.

In the model of Ermentrout-Cowan only position and orientation were considered and symmetry properties were imposed to the facilitation patterns. In \cite{SCS}, for the same features it was proposed to choose as a kernel $\Gamma$ the fundamental solution of a Fokker Planck equation deduced from the Euclidean symmetries. In the next section, two numerical simulations will be performed using equation (\ref{giac}). In the first one, we will model the propagation of boundaries in motion using the connectivity kernel (\ref{eq:densityg0}). In the second experiment we will model instead the propagation of trajectories, using the 5D kernel (\ref{eq:FPg0}).

\section{Numerical simulations}\label{sec:numerical}
In this section some numerical simulations will be performed to test the reliability of the kernels computed in Section \ref{sec:diffusion} when introduced in the activity model (\ref{giac}) and the results will be compared to some psychophysical and physiological experiments.

\subsection{The feedforward and extracellular activity in response to a stimulus}

For the subsequent numerical simulations, we measured in every discrete point $z_j = (x_j, y_j, t_j)$ of a stimulus $f$ of size $n_x \times n_y \times n_t$, its local energy of orientation and speed by convolving the input image sequence with a pre-determined bank of Gabor filters of fixed spatial frequency centered at the points $\eta_i = (x_i,y_i,t_i,\theta_i,v_i)$, thus discretizing Eq. (\ref{eq:filtering}) so to have
\begin{equation}\label{eq:discretefiltering}
\begin{array}{l}
F(\eta_i) = \Big| \displaystyle{\sum_{j=1}^{N}} \, e^{-i( p_1 (x_j-x_i) + p_2 (y_j-y_i) - \nu_i (t_j - t_i))}\\
\hspace{56pt} e^{- \frac{|x_j-x_i|^2 + |y_j-y_i|^2}{\sigma_x^2}-\frac{(t_j-t_i)^2}{\sigma_t^2}} f(z_j) \Big|^2
\end{array}
\end{equation}
where $p_1$ and $p_2$ are the component frequencies $|p| \cos\theta_i$ and $|p| \sin\theta_i$.

We chose the spatio-temporal frequencies and the maximum velocity value $v_{m}$ represented in the Gabor filter set depending on the stimulus, so that given the couple $(|p|, v_m)$ the filters have a maximum temporal frequency of $\nu_m = |p| v_m$. We normalized the filters so that the response could range from $0$ (in regions with no changes in luminance) to $1$ (square plane waves sharing the filter parameters $|p|$ and $\nu_i$), corresponding to a normalization of the stimulus contrast. The Gabor scale parameters will follow the relations
\begin{equation}\label{eq:gabscales}
\left\{
\begin{array}{rcl}
\sigma_x & = & \frac{2.5 \pi }{4 |p|}\vspace{4pt}\\
\sigma_t & = & \frac{\pi}{2 \nu_m}\vspace{4pt}\\
\end{array}
\right.
\end{equation} 
whose meaning is to approximately have $2.5$ spatial subregions under the Gabor Gaussian envelope, and a variable $\nu_i$ that allows the filter with maximum velocity to cover one wavelength within the Gabor active time interval (see \cite{CBS,DeAngelis1993a,Ringach2002} for the physiological justifications).

The thresholded feedforward input $F_{T}$ is then computed following Eq. (\ref{eq:v1thresholding}), representing the mean neural extra-cellular activity in response to the stimulus. In the subsequent analyses, the values to which we set the parameters $\mu$ and $\beta$ in the above stated formulas will be explicitly specified. It is worth noting that these parameters are not due to physiological data, due to the difficulty in finding some reference in the literature, but are chosen as to have a computationally reasonable number of non-null measurements $F_{T}$.

\subsection{Experiment 1 - Contours in motion}

\begin{figure*}
\centering
\includegraphics[width=\textwidth]{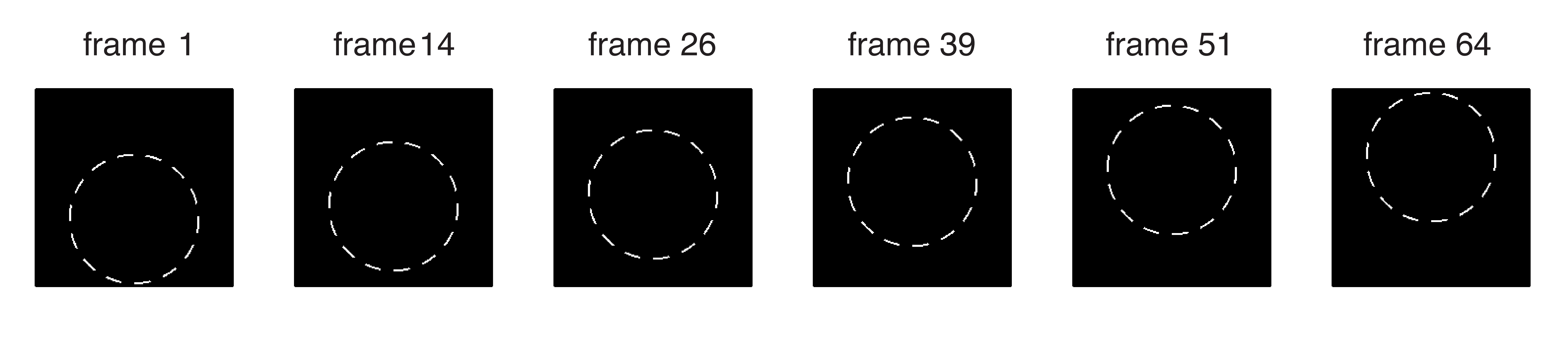}
\caption{The stimulus used in Experiment 1.}\label{fig:rainstim}
\end{figure*}

For the first simulation we use a dashed circle in linear motion moving upwards within the visual space. The image size is $200\times200$ pixels and the image sequence is composed of $64$ frames. The circle segments are approximately $2$ pixels wide and the circle is moving with a uniform speed of $0.5$ pixels per frame (see Fig. \ref{fig:rainstim}). To get a measurement of local spatio-temporal features, we convolve the stimulus with a set of 3-dimensional Gabor functions with (\ref{eq:filtering}). The spatial frequency parameter is set to $|p| = \pi/2$, so that the Gabor moving subregions match the width of the segments. The maximum local velocity represented in the filter set is $v_m = 1$ pixel per frame, while spatial and temporal scale parameters are calculated accordingly with (\ref{eq:gabscales}). After having convolved the stimulus with the Gabors, we model the neural extra-cellular activity by using (\ref{eq:v1thresholding}) with parameters $\mu = 10$, $\beta = 0.5$, obtaining $F_{T}$.

We select the 4-dimensional subset corresponding to the $32$nd frame $F_{T_v}(\xi_i) = F_{T}(x_i,y_i,32,\theta_i,v_i)$, where $\xi_i = (x_i,y_i,\theta_i,v_i)$, thus neglecting behaviors over time. This is because the model given by the stochastic model of connectivity distribution $\Gv$ is stationary, even if it is relying on spatio-temporal information, and is compatible with the dimensionality of $\N$ in (\ref{enne}).

Following the assumptions that the measurements $F_{T_v}(\xi)$ can model the output of a cortical direction- and speed-selective cell (or of a neural population that are selective to the same visual features), the overall continuation probability can be obtained with a discretized version of the cortico-cortical facilitation pattern (\ref{eq:cc}) with the kernel $\Gv$:
\begin{equation}\label{eq:gconv_v}
P(\xi_i) = \sum_{\xi_i' \in \N} \Gv(\xi_i|\xi_i') F_{T_v}(\xi_i')
\end{equation}
where $\Gv(\xi_i|\xi_i')$ is calculated for every fiber vector $(\theta_i, v_i)$ on the same discretized domain $\N$ of the lifted stimulus, using the stochastic approach described in Section \ref{sec:stochasticX2X4X5}. The structure of (\ref{eq:gconv_v}) is that of a discrete version of a group convolution with respect to the composition law (\ref{eq:compositionV}), due to the symmetry (\ref{eq:symmetryV}) of the kernel and the left invariance of the Lebesgue measure on $\N$.

It is worth noting that the parameters $\kappa$ and $\alpha$ governing the diffusion process are related to the maximum spatio-temporal curvature of an illusory contour.
Due to the lack of reference physiological data, the values for the diffusion coefficients were chosen in such a way that at the final value taken by the evolution parameter of the stochastic paths, the mean square displacements of the fiber variables $(\theta,v)$ are $(\pi, v_m/2)$, that means $\kappa = 2$ and $\alpha = 1$.

% Here, we have chosen the values for the diffusion coefficients  manually in order to have the best results, because of the lack of reference physiological data. In particular, the parameters were set so that at the final value taken by the evolution parameter of the stochastic paths, the mean square displacements of the fiber variables $(\theta,v)$ are $(\pi, v_m/2)$.

\begin{figure}
\centering
\subfloat[]{\includegraphics[width=0.48\textwidth]{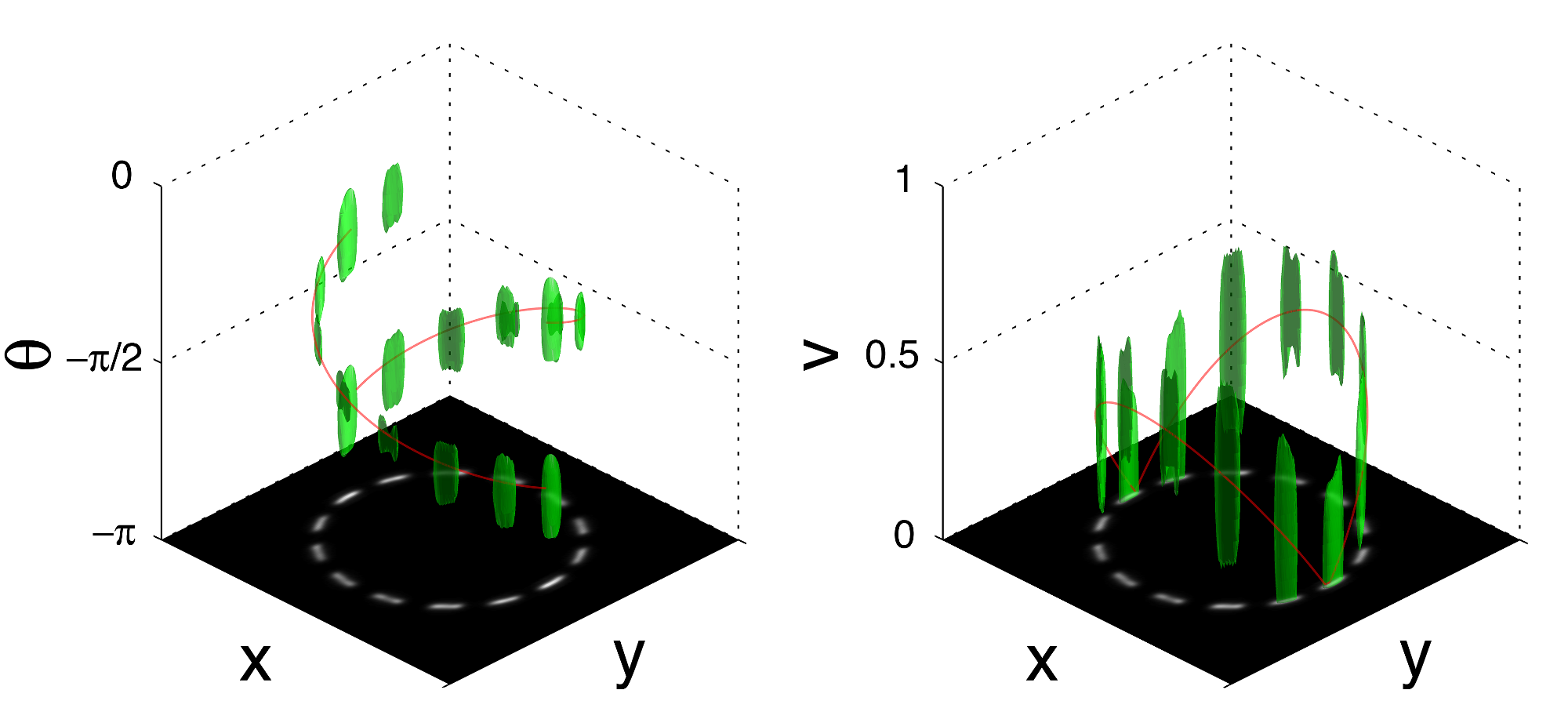}} \quad
\subfloat[]{\includegraphics[width=0.48\textwidth]{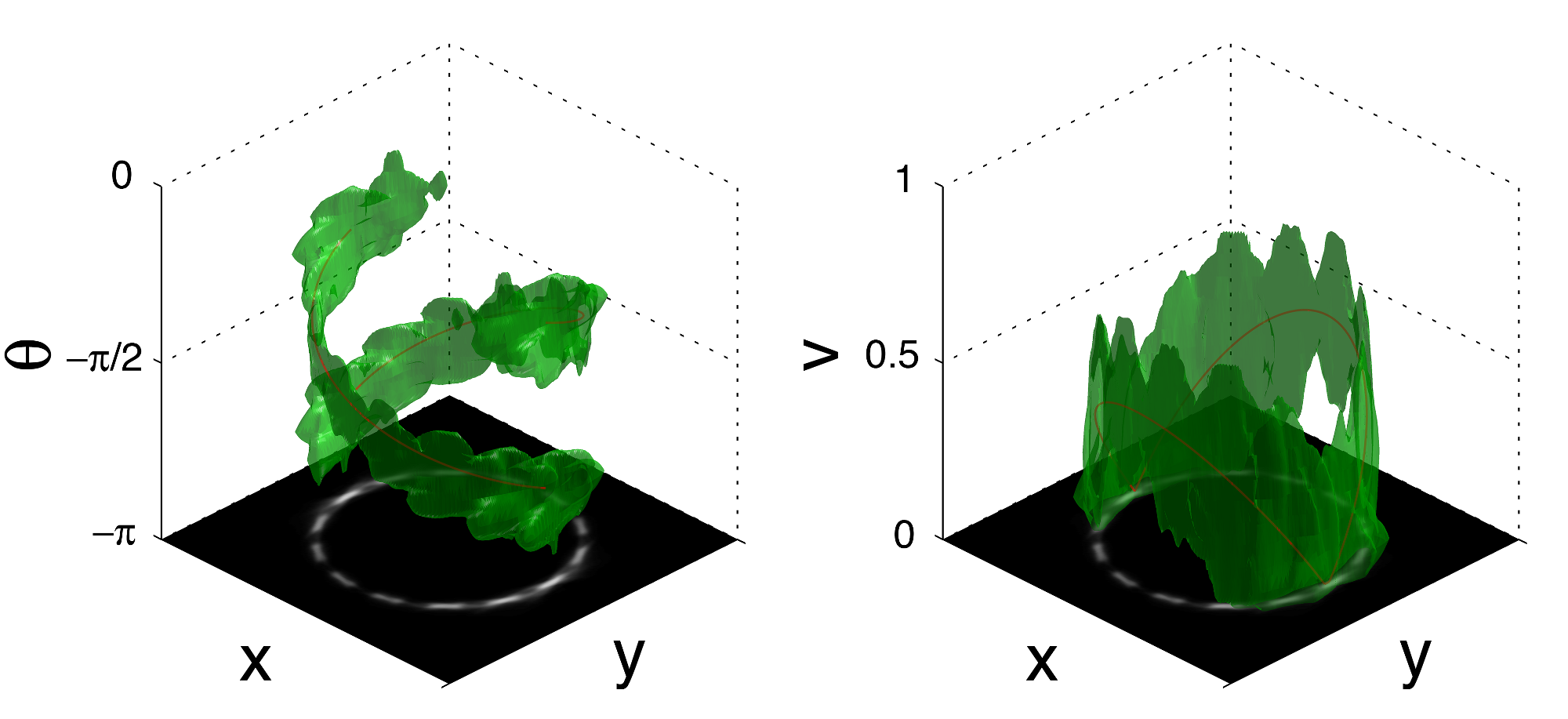}}
\caption{Results of Experiment 1: in the left subplots, measurements of local direction of motion $\theta$; in the right subplots, measurements of local velocity $v$. The theoretical values for a circle translating at a speed of $0.5$ pixels per frame are shown with a red line. a) Isosurface plot of $F_{T}$ (isovalue: $0.2$). The information about local orientation and velocity is correctly retrieved, but activity remains clustered within disjoint regions of the domain. b) The modeled extra-cellular activity $F_{0}$ after horizontal propagation on the submanifold $\N$ (facilitation strength: $c_f = 40$).}\label{fig:rainsym}
\end{figure}

Finally, the population activity is computed by
\begin{equation}\label{eq:weigh_fac}
F_{0}(\xi_i) = S( F(\xi_i) + c_f P(\xi_i))
\end{equation}
that is a discretized version of (\ref{giac}) where $c_f$ is a coefficient governing the total strength of the excitatory connections. In Fig. \ref{fig:rainsym} are shown the results of the simulation. In the top row we have plotted a iso-level surface of the normalized and thresholded output measurements $F_{T_v}$. The process correctly lifts the stimulus around its theoretical values in the domain $\N$. In the bottom row, the iso-level surface of the global cortical activity $F_{0}$ is shown. Note that continuation is performed within the whole manifold $\N$, propagating the local speed and direction cues at the end of the segments, and interpolating the data to define a 4-dimensional set that carries information even in the visual space between the segments, where no changes in contrast were originally present in the stimulus. These results are coherent with recent psychophysiological findings, where it is stated that a global shape in motion is better perceived if local velocity changes smoothly along its contour \cite{Rainville2005}, and that the integration of the motion of a partially occluded object is facilitated when its visible contours define closed configurations \cite{Lorenceau2001}.

\subsection{Experiment 2 - Motion integration}

\begin{figure*}
\centering
\includegraphics[width=\textwidth]{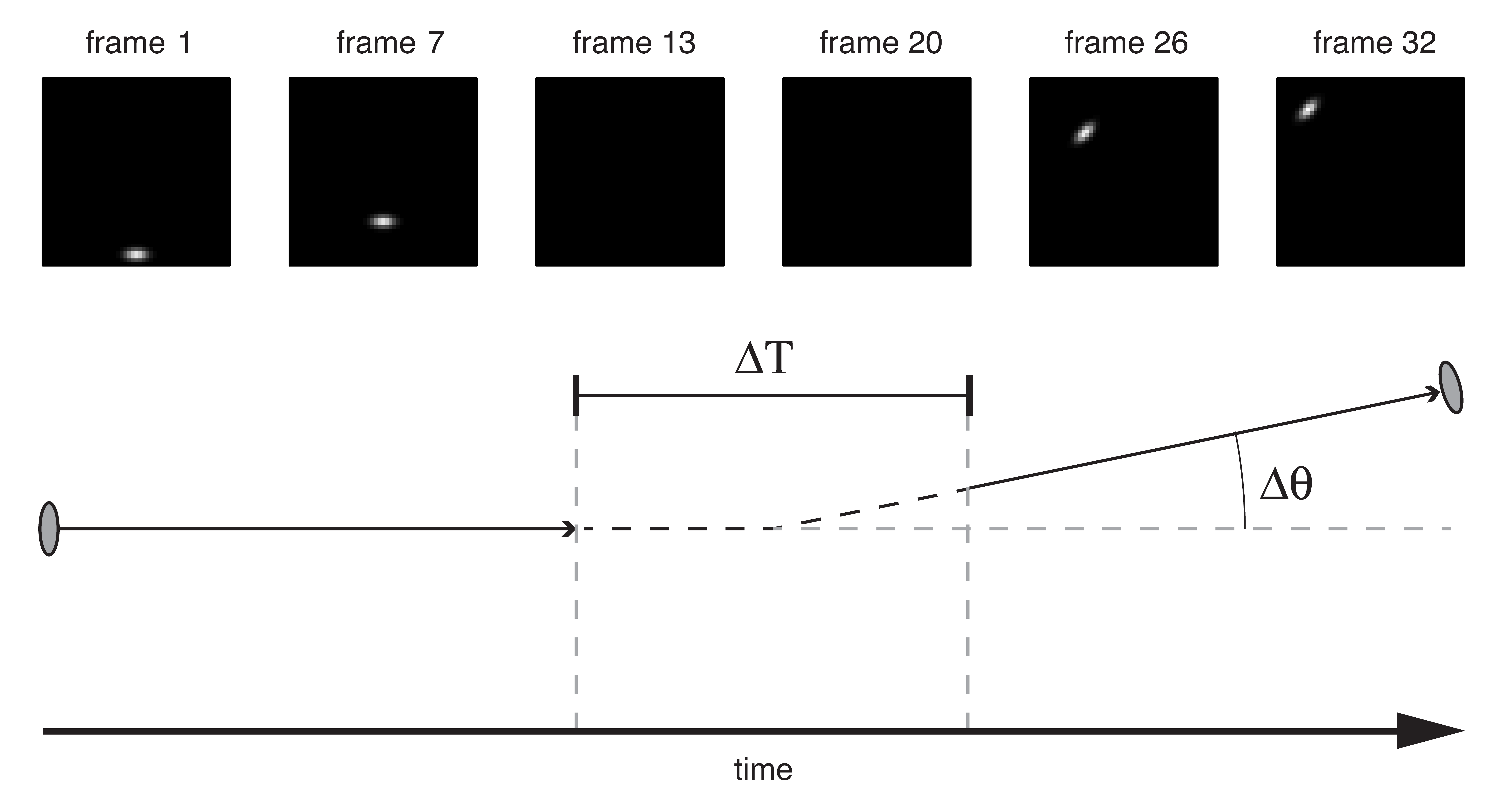}
\caption{An example stimulus sequence, and the paradigm used to generate the stimuli in Experiment 2.}\label{fig:trajstim}
\end{figure*}

The stimuli that we will use in this simulation will be several instances of an object in motion along a certain direction, that disappears at a given time position $t_1$, and reappears at $t_2$ with a direction of motion changed by $\Delta \theta$ in a position that is coherent to a piecewise continuous trajectory (Fig. \ref{fig:trajstim}). It is a well documented fact that humans tend to perceive the two trajectories as pertaining to a single unit just for small values of $\Delta T = t_2 - t_1$ and $\Delta \theta$, an effect that is commonly referred to as motion integration. In particular, we know from many psychophysiological experiments that the chances of detecting a straight or curvilinear trajectory in noise increase with stimulus duration and with direction coherence \cite{Watamaniuk1995, Verghese1999}.

To do so, we created multiple instances from the stimulus paradigm described above, and shown in Fig. \ref{fig:trajstim}, assigning different values to the parameters $\Delta T$ and $\Delta \theta$. It is worth noting that we are using an elongated object whose axis of motion is always orthogonal to its major axis of eccentricity. The choice of this particular stimulus is naturally explained by the geometry, as the points $\in \M$ that are connected by its continuation property implicitly define a local orientation and a direction of motion that is orthogonal to that orientation. %Furthermore, an elongated object will give rise to neater responses after the convolution with our filter set.
The Gabor filtering of a different kind of stimulus (for example a moving dot) would indeed measure high responses also for those velocity projections that are coherent with the real axis of motion, due to the motion streak effect \cite{Geisler1999}, requiring a different neural model to detect unambiguously the direction of motion. Thus, even if the majority of the psychophysiological experiments use moving dots as stimuli, the validity of the following experiment is not influenced, as exploring the exact physiological implementation of motion integration is out of the scope of this paper.

The image sequences to process are $51\times51$ pixels wide and are composed of $102$ frames. The value of eccentricity of the moving ellipsoidal object is set to $2$, and its velocity is set to $0.5$ pixels per frame. Each stimulus instance is uniquely identified by the couple $(\Delta T, \Delta \theta)$. We convolve each sequence with a set of Gabor filters with $v_m=1$ and $p$ set to match width of the object's minor axis, and again, we threshold the output according to (\ref{eq:v1thresholding}) using parameters $\mu = 20$, $\beta = 0.7$, obtaining the set of measurements $F_{T}$. In the previous simulation we have taken a temporal slice of this output, and then we have propagated the activity using the connectivity kernel $\Gv$. Similarly, here we will propagate the activity present in $F_{T}$, but without discarding time. The corresponding facilitation pattern $P(\eta_i)$ is then obtained convolving it with the 5-dimensional kernel $\Go$ defined on the same discrete domain $\M$ of the lifted stimulus using the same approach that led to (\ref{eq:gconv_v})
\begin{displaymath}
P(\eta_i) = \sum_{\eta_i' \in \M} \Go(\eta_i|\eta_i') F_{T}(\eta_i')
\end{displaymath}
and the total population activity  $F_{0}(\eta_i)$ is computed following (\ref{eq:weigh_fac}) with $c_f = 20$. Since this is a discretization of (\ref{eq:cc}), again we are considering a convolution-type operator.
\begin{figure*}
\centering
\subfloat[]{\includegraphics[width=.48\textwidth]{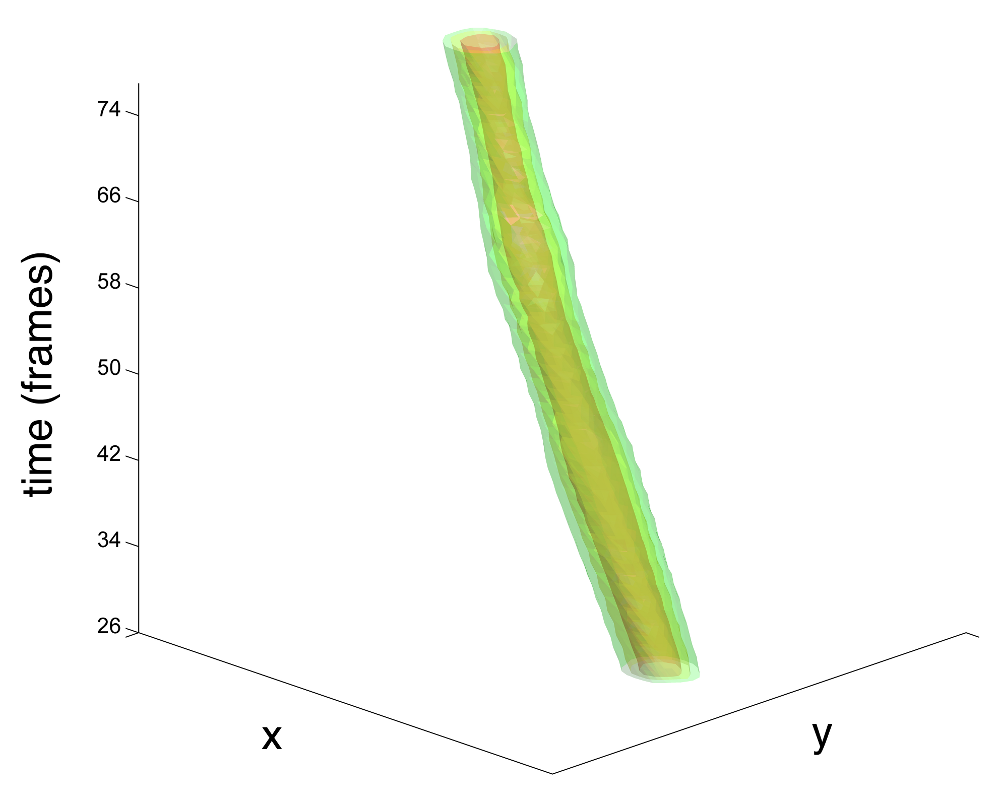}} \quad
\subfloat[]{\includegraphics[width=.48\textwidth]{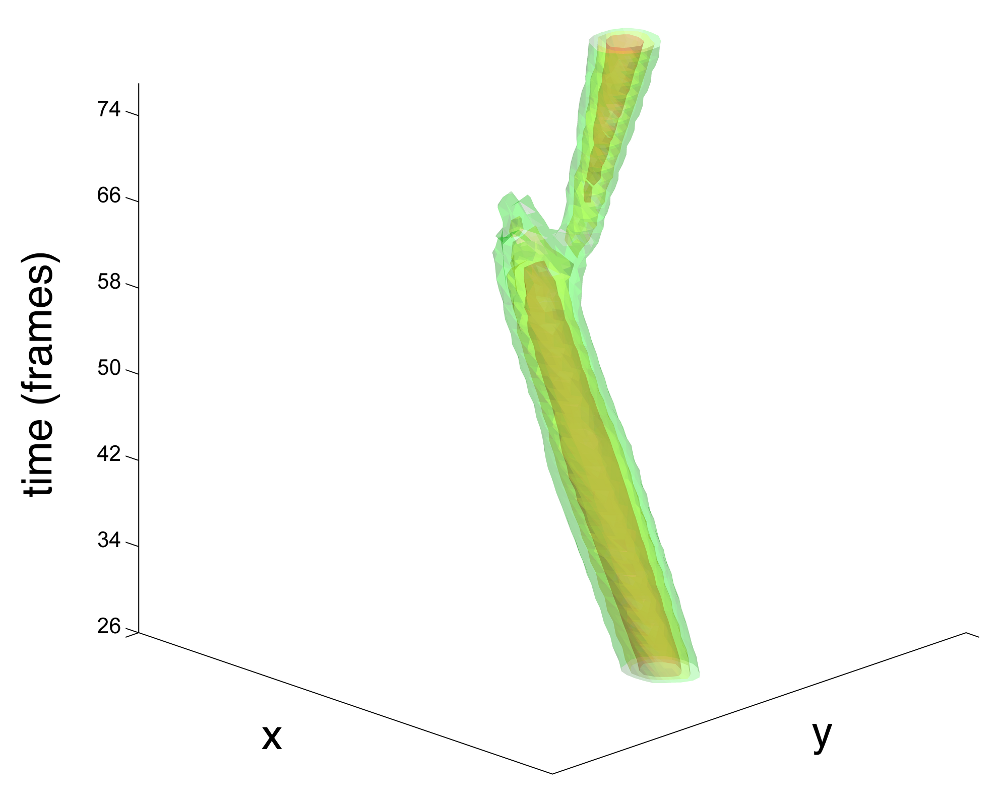}}
\caption{A comparison between the spatio-temporal continuation capabilities of the model, applied to stimulus instances with different values of $\Delta \theta$. The plot shows three isosurfaces of $\max_{(\theta,v)} F_{0}$ (isovalues 0.9, 0.5, 0.1 in red, yellow an green): a) $\Delta T = 12$, $\Delta \theta = \frac{\pi}{6}$, b) $\Delta T = 12$, $\Delta \theta = \frac{5\pi}{12}$. The probabilistic continuation given by the geometry in the configuration with the smallest change in direction of motion allows a smooth trajectory completion. The same is not true for the stimulus with an higher value of $\Delta \theta$: even if a weak response connecting and interpolating the two parts of the stimulus is still present, the strongest facilitation component do not deviate from the early straight path of the object.}\label{fig:trajsym_xyt}
\end{figure*}

\begin{figure*}
\centering
\subfloat[]{\includegraphics[width=.48\textwidth]{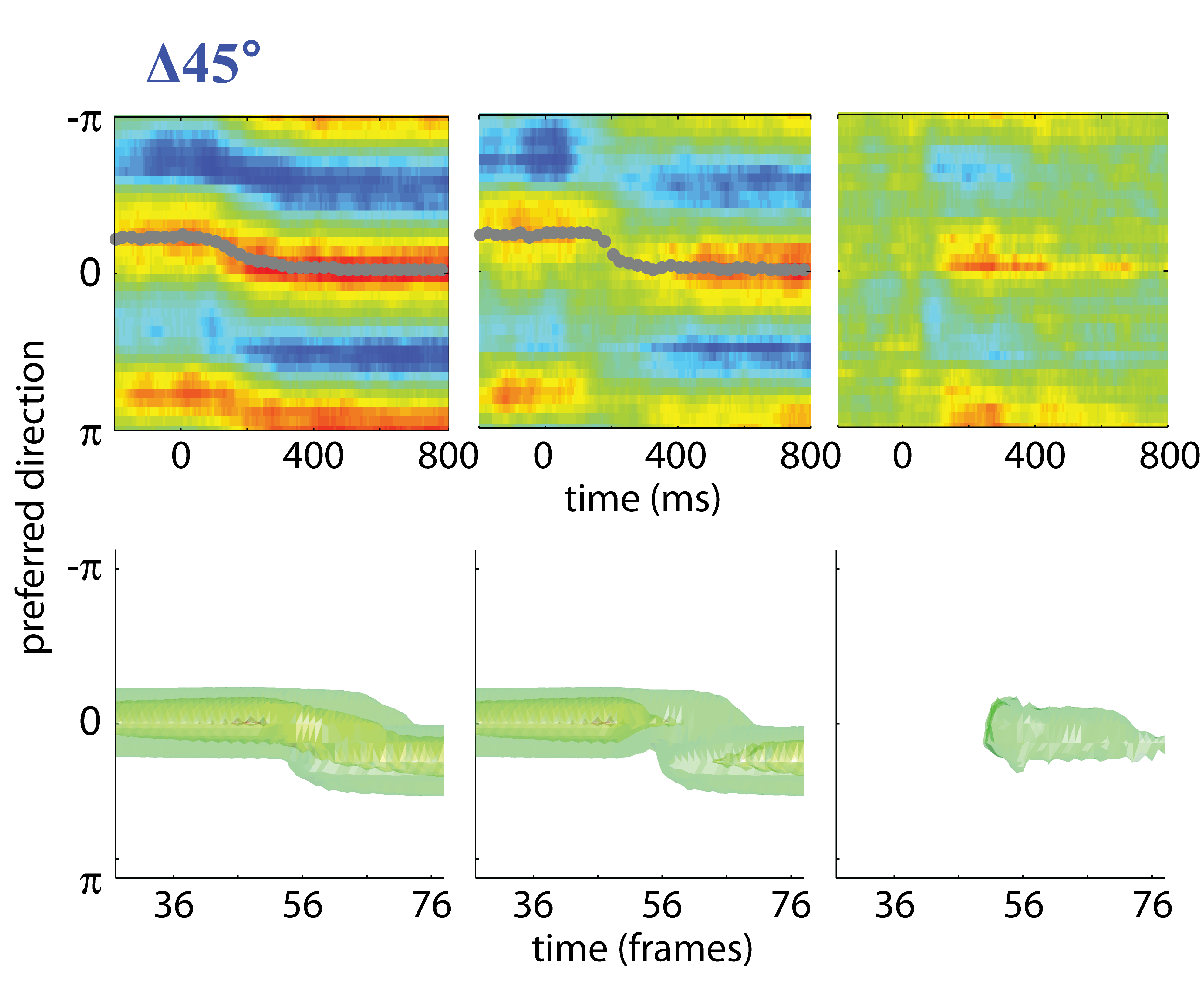}} \quad
\subfloat[]{\includegraphics[width=.48\textwidth]{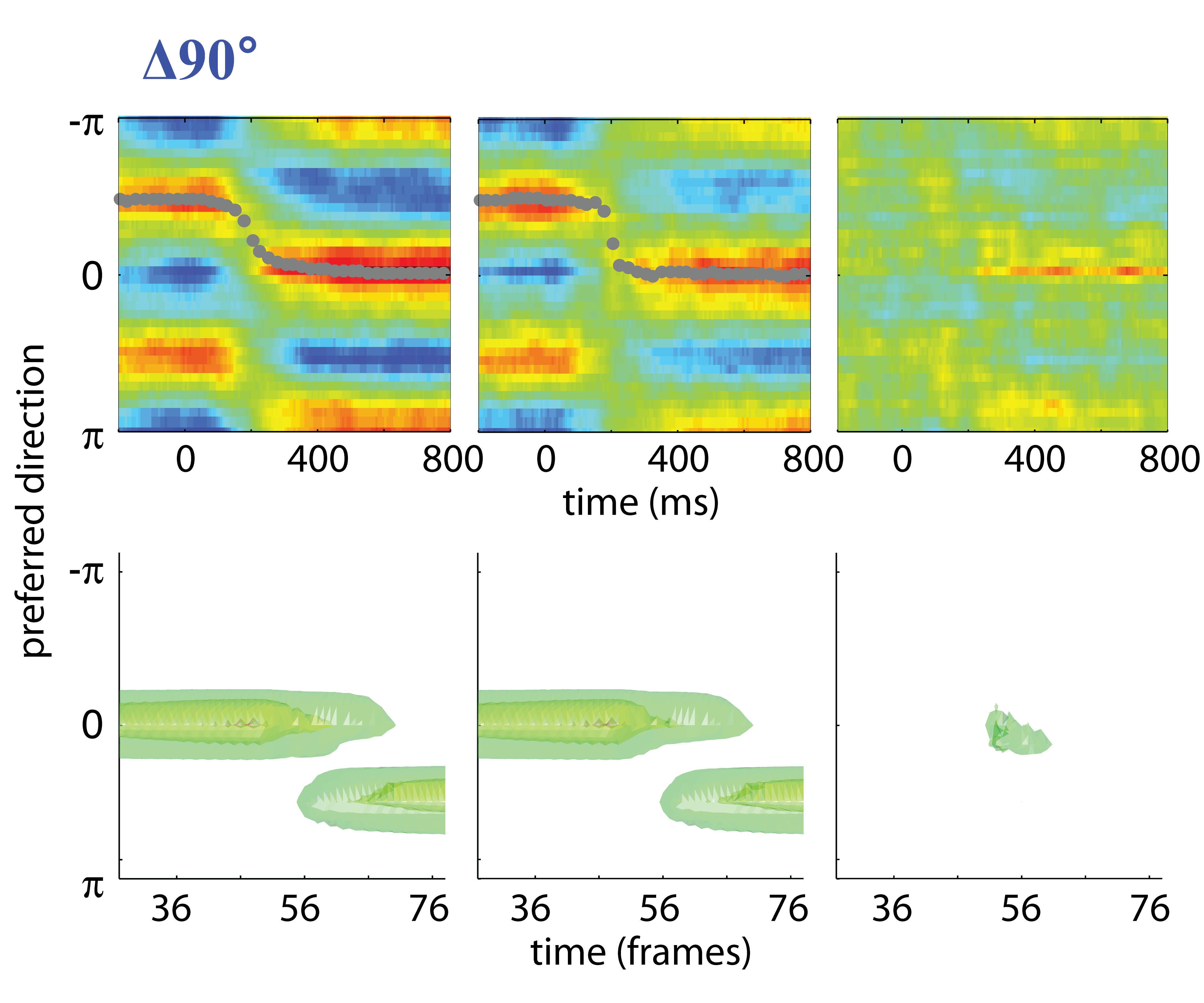}}
\caption{A direct comparison between the electrophysiological recordings made in \cite{Wu2011} and the results obtained with our model for the same stimulus parameters: in both cases $\Delta T = 0$, while $\Delta \theta$ is equal to $\frac{\pi}{4}$ (a) and $\frac{\pi}{2}$ (b). To compare responses, we integrate over $x$ and $y$ and depict three isosurfaces with isovalues taken as 0.9, 0.5, 0.1 (in red, yellow and green respectively) times the maximum value among the outputs shown. Note that, disregarding the weaker physiological responses to opposite directions, the non-linear facilitating behavior is very well reproduced by our model.}\label{fig:trajsym_fac}
\end{figure*}

The results obtained by processing two stimulus configurations are shown in Fig. \ref{fig:trajsym_xyt}, where only the central portion of the full temporal domain is plotted, to highlight the different effects that activity propagation has in the two cases. When the change in direction of motion of the object is under a certain threshold, the trajectory is completed smoothly, yielding a strong activation in the temporal interval when the object is not present in the scene. If the value of $\Delta \theta$ increases, however, the spatio-temporal curvature of the optimal subjective trajectory progressively becomes too high to be completed by the connectivity, even if a weak response linking the two parts of the stimulus could yet be present.

It is worth noting that due to the thresholding stage, this subjective interpolation effect is strongly non-linear, as it cannot be explained by the sum of the responses to the first ($t_i<t_1$) and second ($t_i>t_2$) part of the stimulus alone (see Fig. \ref{fig:trajsym_fac}). This is coherent with the findings in the work of Wu et al, where electrophysiological experiments recorded a similar non-linear behavior of trajectory interpolation for small values of $\Delta \theta$ \cite{Wu2011}. In that work, the cause of this effect is left unexplained. Even if the stimulus paradigms are slightly different (they use a field of moving dots with a null value of $\Delta T$), the result of our simulations allow us to make some speculations. 

\begin{figure}
\centering
\includegraphics[width=.65\textwidth]{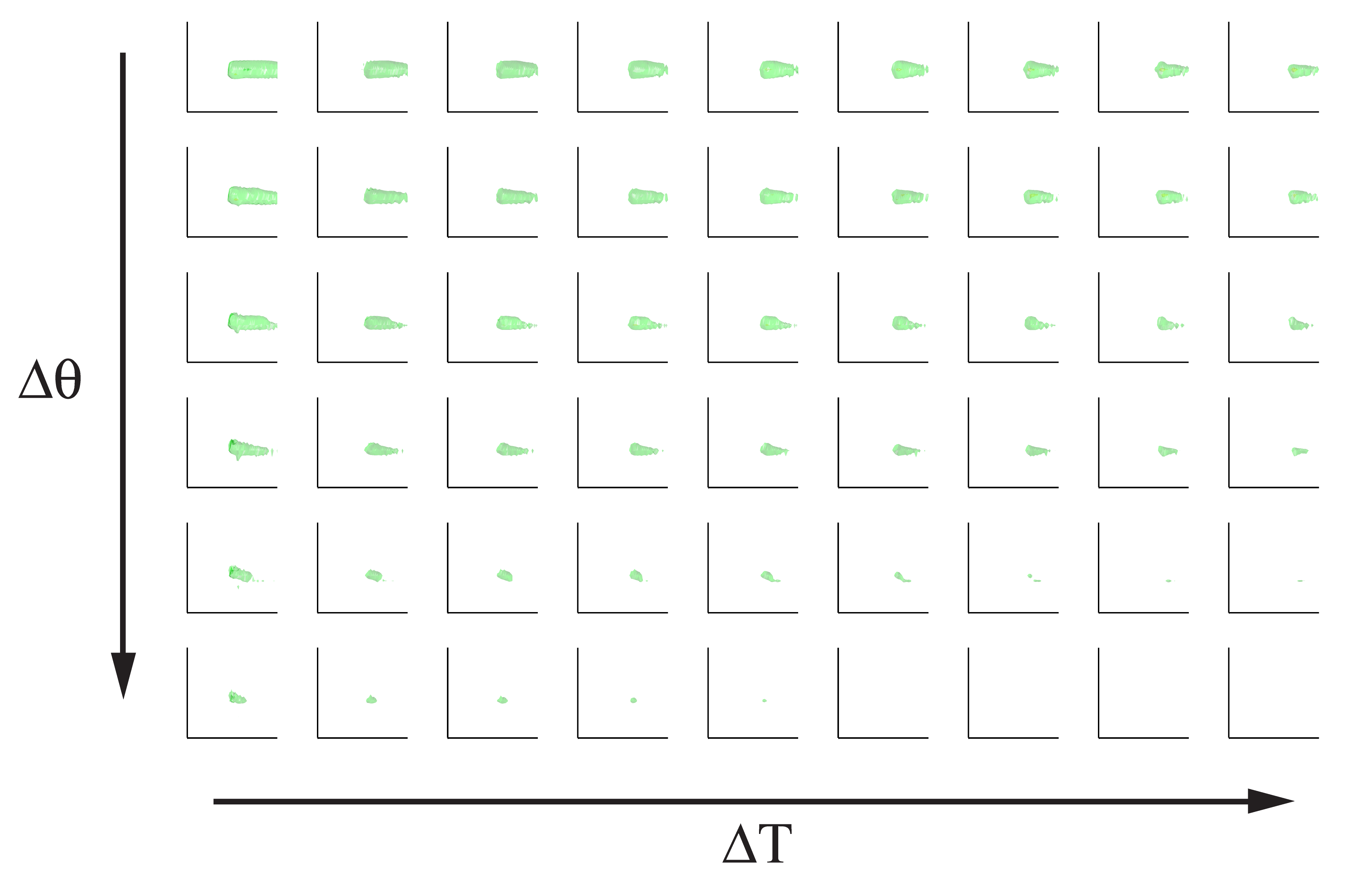}
\caption{The non-linear effect of the facilitation by varying the stimulus parameters $\Delta T$ and $\Delta \theta$. Each subplot shows the isosurface of $F_{fac}$ relative to different stimulus instances (isolevels calculated as in Fig. \ref{fig:trajsym_fac}). The non-linear effect of the facilitating kernel $\Go$ remains appreciable just for small values of $\Delta T$ and $\Delta \theta$, rapidly decaying as the value of the characteristic stimulus parameters increases.}\label{fig:trajsym_fac2}
\end{figure}

A good qualitative description of the effects that the modeled excitatory connectivity has on stimulus response can be found in Fig. \ref{fig:trajsym_fac2}, where for every stimulus configuration we plot the difference
\begin{equation}\label{eq:diff_fac}
F_{fac} = F_{0}(S_3) - F_{0}(S_1) - F_{0}(S_2)
\end{equation}
where $S_3$ is the full stimulus and $S_1$ and $S_2$ represent the first ($S1(t_i>t_1)=0$) and the second ($S2(t_i<t_2)=0$) part of the stimulus. The visualization of the output $F_{fac}$ highlights the role that a trajectory-specialized cortical connectivity could have in performing tasks of motion integration. The spatio-temporal extension of $F_{fac}$ over $\M$ gets smaller for higher values of $\Delta T$ and $\Delta \theta$. The decay of the facilitation effect is coherent with the observations made in \cite{Watamaniuk2005} (where $\Delta \theta = 0$) and \cite{Wu2011} (where $\Delta T = 0$), even if the experiments are methodologically different. Moreover, even if, coherently with our results, some psychophysiological experiments showed that a broken trajectory in noise is easily detectable \cite{Watamaniuk1995n, Scholl1999}, as far as we know little has been done to explore the effect of changing the duration $\Delta T$ of the temporal occluding gap.

\subsection{Discussion}

Regarding the implementation of our model, we would like to highlight functional meaning of the parameter couple $(\kappa, \alpha)$, driving the diffusion along the fiber variables $(\theta, v)$ when calculating the kernels $\Gv$ and $\Go$. In particular, the parameter $\kappa$ seems to be strongly related to the maximum perceived curvature of illusory contours and trajectories, while $\alpha$ sets the maximum rate of change of local velocity along admissible subjective contours (trajectories). It would be interesting to try to fine-tune the parameters of the model in order to reproduce quantitatively as precisely as possible the psychophysiological findings found in literature. We aim to carry out this analysis in a future paper.

It is worth noting that the physiological correlates of the first simulation are not well documented and, as far as we know, no one has yet studied the neural activity in cortical regions responding to subjective contours in motion. We suppose that it would be a difficult issue to address, as multiple cortical layers, as V1, V4 and MT/V5, may be involved \cite{Lee2001}. Recently, some significative results have been obtained with static illusory contours, using promising electrophysiological techniques \cite{Pan2012}. Those kinds of experiments, targeted to the detection of illusory motion contours, could give additional clues about the neural computation that governs the influence of spatio-temporal features in the detection of moving shapes and boundaries.

As discussed in the beginning of Section \ref{sec:truelife}, a possible physiological implementation of the geometry used in the second numerical simulation, could be a trajectory-specialized cortico-cortical connectivity between neurons in the higher visual areas. The advantages of having such an excitatory connectivity implemented in the visual cortex could be very important in performing complex cognitive tasks. For example, we assume that the brain could use the improvements in contour detection with respect to the only $SE(2)$ related functional architecture, due to an additional information, that of velocity, that indeed in practical situation is a feature that is coherent on objects. The influence that the connectivity arising from our geometrical model has in the processing of visual tasks such as spatio-temporal grouping or segmentation will be the subject of a future paper.

\section{Conclusions}

In this paper we proposed a model of cortical functional architecture for the processing of spatio-temporal visual information.  Motion features are detected first by simple and complex cells with RP modelled by 2D+time Gabor filters. Linear filtering with such a profiles lifts the visual stimulus from the base space $\R^3$ to the phase space $\R^6$ comprising spatial and temporal frequencies, in which a Liouville form can be defined. 

We defined a 5D phase space with fixed frequency by taking a reduction of the previous differential form. Then, exploiting the commutation properties of the horizontal basis, we regarded the tangent space of the contact manifold $\M$ as a possible constraint acting on the connectivity between points, giving the definition of admissible integral, or horizontal, curves.
Possible linear combinations of the horizontal basis that are compatible with the definition of admissible curve have been studied, in order to model the possible lifting of the visual stimuli as association fields, in the sense of \cite{Field1993}. We considered the corresponding deterministic integral curves for two modeling limit cases: contours in motion and trajectories of a point in motion. 

Then we considered horizontal stochastic paths, i.e. trajectories of points that always move along the tangent space of $\M$ and are allowed to change the value associated to the fiber variables $(\theta,v)$ in a random, equidistributed way. We have seen that the resulting integrations over the evolution parameter of the paths coincide with the kernels of the Fokker Planck operators defined on the geometries. Both deterministic curves and stochastic kernels inherit the symmetries of the non associative reduced Galilean group as described in Appendix.

After a  discussion about the compatibility of our theoretical framework with various phenomenological and psychophysiological findings on visual perception and cognition reported in the literature, we introduced the stochastic kernels as facilitation inducers in a neural population activity model.

Suitable numerical simulations have been carried out, by processing pre-determined artificial stimuli with the neural population activity model previously descri\-bed, showing that the Fokker Planck stochastic kernels endow the model with the capability of completion and continuation of contours in motion and trajectories of points, coherently with the phenomenological experiments of Rainville \cite{Rainville2005} and Wu et al \cite{Wu2011}. \\
In conclusion, results have shown that the proposed functional geometry is compatible with existing  psychophysical and physiological experiments, even if a complete knowledge about the effective neural implementation needs supplementary empirical data.

\section*{Appendix: Relations with the Galileian group}\label{sec:Galilei}

The geometry of $(\M,\omega)$ can be related with the classical Galilei group of motions on the plane, when choosing an appropriate local frame of reference.

The two dimensional Galilei group is the semidirect product Lie group $\G = (\R^2_q \times \R_s) \ltimes (\R^2_u \ltimes S^1_\theta)$ defined by the composition law \cite{Sorba}
\begin{equation}\label{eq:galileilaw}
\begin{array}{l}
(\q,s,\theta,u) \cdot (q',s',\theta',u')\\ = (R_\theta q' + u s' + \q, s' + s,\theta' + \theta, R_\theta u' + u)
\end{array}
\end{equation}
where $(\q,s)$ provides space-time translations, $R_\theta$ is a counterclockwise planar rotation of an angle $\theta$ and $u = \textstyle{\binom{u_1}{u_2}}$ is the velocity vector responsible of the purely Galilean transformation, or boost, expressed in Cartesian coordinates.

The manifold $\M$ can be identified with a subset of this group by identifying the direction of velocity with the angle of the rotation. Setting the velocity in polar coordinates
\begin{displaymath}
u = \left\{
\begin{array}{rcl}
u_1 & = & v \cos\varphi\\
u_2 & = & v \sin\varphi
\end{array}
\right.
\end{displaymath}

then we can identify $\M$ with the set of points
\begin{displaymath}
\M = \{(\q,s,\theta,R_\theta \mathfrak{v})\}
\end{displaymath}
that is, using the identification $\theta = \varphi$. Since the angle $\theta$ is interpreted as the local direction of a boundary, then this amounts to select the apparent velocity orthogonal to it.

This manifold is not invariant under the Galilei composition, but it is invariant with respect to the smooth composition law (\ref{eq:composition}). This composition law is not associative, but possesses a trivial neutral element $e = 0$ and each $\eta \in \M$ possesses a left inverse element, i.e. $\eta^{-1_L} \in \M$ such that $\eta^{-1_L} \odot \eta = 0$
\begin{displaymath}
(\q,s,\theta,v)^{-1_L} = (-R_{-\theta} (\q - \textstyle{\binom{v}{0}} s), -s,-\theta,-v)
\end{displaymath}
and for these reasons we can consider it a reduced Galilei non associative group \cite{Sabinin}.

In particular, standard invariance arguments apply to the solution to (\ref{eq:FPg0}), that possesses the symmetry
\begin{displaymath}
\Go(\eta|\eta_0) = \Go(\eta_0^{-1_L}\odot \eta|0)\quad \forall \ \eta,\eta_0 \in \M\ .
\end{displaymath}
Analogously, the fan $\Sigma(\eta_0)$ of integral curves (\ref{eq:associationT}) with constant coefficients starting from a point $\eta_0$ can be obtained from the one starting from $0$ by
\begin{displaymath}
\Sigma(\eta_0) = \eta_0^{-1_L}\odot \Sigma(0)\ .
\end{displaymath}

\end{document}